\documentstyle[bezier]{article}

\makeatletter
\def\eqnarray{%
  \stepcounter{equation}%
  \let\@currentlabel=\theequation
  \global\@eqnswtrue
  \global\@eqcnt\z@
  \tabskip\@centering
  \let\\=\@eqncr
  $$\halign to \displaywidth\bgroup\@eqnsel\hskip\@centering
  $\displaystyle\tabskip\z@{##}$&\global\@eqcnt\@ne
  \hfil$\displaystyle{{}##{}}$\hfil
  &\global\@eqcnt\tw@$\displaystyle\tabskip\z@{##}$\hfil
  \tabskip\@centering&\llap{##}\tabskip\z@\cr}
\makeatother

\makeatletter
  \renewcommand{\theequation}{%
        \thesection.\arabic{equation}}
  \@addtoreset{equation}{section}
\makeatother

\begin{document}

\newtheorem{th}{Donotwrite}[section]

\newtheorem{theorem}[th]{Theorem}
\newtheorem{proposition}[th]{Proposition}
\newtheorem{conjecture}[th]{Conjecture}
\newtheorem{lemma}[th]{Lemma}
\newtheorem{corollary}[th]{Corollary}
\newtheorem{remark}[th]{Remark}
\newtheorem{example}[th]{Example}

\newfont{\germ}{eufm10}
\newfont{\slsmall}{cmsl8}

\def\B{{\cal B}}
\def\bbar{\overline{b}}
\def\btilde{\tilde{b}}
\def\cd{\cdots}
\def\eps{\epsilon}
\def\et#1{\tilde{e}_{#1}}
\def\ft#1{\tilde{f}_{#1}}
\def\geh{\goth{g}}
\def\goth#1{\mbox{\germ #1}}
\def\L{{\cal L}}
\def\La{\Lambda}
\def\la{\lambda}
\def\ol#1{\overline{#1}}
\def\ot{\otimes}
\def\P{{\cal P}}
\def\pbar{\overline{p}}
\def\Pcl{P_{cl}}
\def\Pcll{(P_{cl}^+)_l}
\def\Proof{\noindent{\sl Proof.}\quad}
\def\qed{~\rule{1mm}{2.5mm}}
\def\slchap{\widehat{\goth{sl}}_{\,n}}
\def\veps{\varepsilon}
\def\vphi{\varphi}
\def\wt{\mbox{\sl wt}\,}
\def\wts{\mbox{\slsmall wt}\,}
\def\Z{{\bf Z}}
\def\Zn{\Z_{\ge0}}

\title{ Character Formulae of $\widehat{sl}_n$-Modules \\
        and Inhomogeneous Paths }

\author{
Goro Hatayama\thanks{
Institute of Physics, University of Tokyo, Komaba, Tokyo 153-8902, Japan},
Anatol N. Kirillov\thanks{
CRM, Universit\'e de Montr\'eal, C.P.6128-succursale Centre-Ville,
Montr\'eal, QC H3C 3J7, Canada \&
Steklov Mathematical Institute, Fontanka 27,
St. Petersburg, 191011, Russia},
Atsuo Kuniba,$\hspace{-1.2mm}^*$ \\
Masato Okado\thanks{
Department of Mathematical Science, Faculty of Engineering Science,
Osaka University, Toyonaka, Osaka 560-8531, Japan},
Taichiro Takagi\thanks{
Department of Mathematics and Physics, National Defense Academy,
Yokosuka 239-8686, Japan}
and Yasuhiko Yamada\thanks{
Department of Mathematics, Faculty of Science,
Kobe University, Rokko, Kobe 657-8501, Japan}
}

\date{}
\maketitle

\begin{abstract}
\noindent
Let $B_{(l)}$ be the perfect crystal for
the $l$-symmetric tensor representation of
the quantum affine algebra $U'_q(\slchap)$.
For a partition
$\mu = (\mu_1, \ldots, \mu_m)$, elements of the tensor product
$B_{(\mu_1)}\otimes \cdots \otimes B_{(\mu_m)}$
can be regarded as inhomogeneous paths.
We establish a bijection between
a certain large $\mu$ limit of this crystal
and the crystal of an (generally reducible)
integrable $U_q(\slchap)$-module, which
forms a large family depending on the inhomogeneity of
$\mu$ kept in the limit.
For the associated one dimensional sums,
relations with the Kostka-Foulkes polynomials are clarified,
and new fermionic formulae are presented.
By combining their limits with the bijection,
we prove or conjecture
several formulae for the string functions,
branching functions, coset branching functions
and spinon character formula of both vertex and RSOS types.
\end{abstract}

{\em Contents}.
{\em 0 Introduction}.  {\em 1 Crystals}.
1.1 Preliminaries.
1.2 Nakayashiki-Yamada's energy function.
1.3 Evaluation of the energy.
{\em 2 Tensor product of crystals}.
2.1 Decomposition of $B(\la)\otimes B_{(l)}$.
2.2 Limit.
2.3 Proof.
2.4 Affine weight.
{\em 3 One dimensional sums}.
3.1 Unrestricted, classically restricted and restricted paths.
3.2 1dsums and Kostka-Foulkes polynomials.
3.3 Limit of 1dsums.
{\em 4 Fermionic formulae}.
4.1 Fermionic formulae of 1dsums.
4.2 Fermionic formulae of the limits.
{\em 5 Discussion}.
5.1 Fermionic string function for arbitrary $X^{(1)}_n$.
5.2 Fermionic form of $X^{(l)'}_\eta(\la)$ for $\la$ non vacuum type.

\setcounter{section}{-1}

\section{Introduction}

Probably, the Kostka-Foulkes polynomial ranks among the most important
polynomials in combinatorics and representation theory. Let $\la,\mu$
be partitions with the same number of nodes. The Kostka-Foulkes
polynomial $K_{\la\mu}(q)$ is defined as the transition matrix which
expresses the Schur function $s_\la(x)$ in terms of the Hall-Littlewood
polynomials $P_\mu(x;q)$: $s_\la(x)=\sum_\mu K_{\la\mu}(q)P_\mu(x;q)$.
(See \cite{Ma} for details.)

Let us consider the affine Lie algebra $\slchap$. We denote by $V(l\La_0)$
the integrable highest weight $\slchap$-module with highest weight
$l\La_0$. Let $\la$ be a partition whose depth is less than or equal to
$n$. We further assume $|\la|$ ($=$ the number of nodes in $\la$) is
divisible by $n$. $\la$ can also be viewed as a level $l$ integral weight
by $(l+\la_n-\la_1)\La_0+(\la_1-\la_2)\La_1+\cd+(\la_{n-1}-\la_n)\La_{n-1}$.
In \cite{Ki}, A.N. Kirillov conjectured the following identity.
\begin{eqnarray*}
\lim_{N\to\infty}&&\hspace{-.4cm}
q^{-E_N}K_{((lN-|\la|/n)^n) + \la,(l^{nN})}(q)
=\sum_j(\dim M^{l\La_0}_{\la-j\delta})q^j,\\
M^{l\La_0}_\mu&=&\{v\in V(l\La_0)\mid e_iv=0\,(i\ne0),\wt v=\mu\}.
\end{eqnarray*}
Here $E_N$ is a known constant. For the definition of $(k^n) +  \la$, see
the beginning of Section \ref{subsec:limit}. It had not been long before
Nakayashiki and Yamada \cite{NY} solved this conjecture. Their idea was
to relate Lascoux-Sch\"utzenberger's charge of a tableau with the so-called
{\em energy} of a path. Once this correspondence is established,
the conjecture is found to be a corollary of the theory of perfect crystals
\cite{KMN1,KMN2}.

The purpose of this paper is to extend their
result to more general setting and elucidate an interplay
among the theory of crystals, the Kostka-Foulkes polynomials,
one dimensional sums, their fermionic formulae
and affine Lie algebra characters.
In a sense this is a far reaching application of the corner transfer matrix
method \cite{ABF} and the Bethe ansatz \cite{Bethe}
in solvable lattice models
where many important ideas came from.
Let us give below an overview of the main contents and results.

In Section \ref{sect1} we recall the definition of the
energy in crystal base theory based on \cite{NY}.

In Section \ref{sect2} we prove our first main Theorem \ref{th:main_th},
which establishes a bijection of
crystals related to a large family of (generally reducible)
integrable highest weight $U_q(\slchap)$-modules.
To explain it more concretely, let us introduce some notations.
Let $B(\la)$ be
the crystal base of the integrable highest weight
$U_q(\slchap)$-module with highest weight $\la$. The symmetric tensor
representation of $U'_q(\slchap)$ of degree $l$ also has a crystal base,
which is denoted by $B_{(l)}$. Let $\mu$ be a partition with signature
$(\mu_1,\mu_2,\cd,\mu_m)$. We consider the tensor product
\begin{equation} \label{eq:tensor_prod}
B_{(\mu_1)}\ot B_{(\mu_2)}\ot\cd\ot B_{(\mu_m)}.
\end{equation}
An element of this tensor product is to be called {\em inhomogeneous}
path, since the degrees $\mu_i$ are not necessarily equal. From the
viewpoint of crystal base theory, Kirillov's conjecture corresponds to
the fact that if $\mu$ is of shape $(l^{nN})$, the crystal
(\ref{eq:tensor_prod}) in the limit $N\to\infty$ is bijective to
$B(l\La_0)$. In contrast to this, we consider in this paper the case
when $\mu$ has the form:
%************************** mu *********************************:
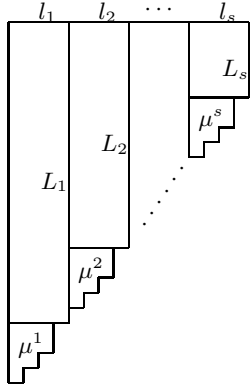
\begin{center}
\unitlength=2mm
{\small
\begin{picture}(16,28)(0,1)
\put(0,25){\line(1,0){16}}
\put(0,25){\line(0,-1){24}}
\put(4,25){\line(0,-1){20}}
\put(8,25){\line(0,-1){15}}
\put(12,25){\line(0,-1){9}}
\put(16,25){\line(0,-1){5}}
\put(2,25.2){\makebox(0,0)[lb]{$l_1$}}
\put(6,25.2){\makebox(0,0)[lb]{$l_2$}}
\put(10,25.5){\makebox(0,0)[b]{$\cdots$}}
\put(14,25.2){\makebox(0,0)[lb]{$l_s$}}
\put(4,15){\makebox(0,0)[rt]{$L_1$}}
\put(8,17.5){\makebox(0,0)[rt]{$L_2$}}
\put(16,22.5){\makebox(0,0)[rt]{$L_s$}}
\put(1.5,3.5){\makebox(0,0){$\mu^1$}}
\put(5.5,8.5){\makebox(0,0){$\mu^2$}}
\put(13.5,18.5){\makebox(0,0){$\mu^s$}}
\multiput(0,1)(4,5){2}{                      %\mu^1,\mu^2
	\multiput(0,0)(1,1){3}{\line(1,0){1}}
	\multiput(1,0)(1,1){2}{\line(0,1){1}}
	\put(3,2){\line(0,1){2}}
	\put(0,4){\line(1,0){4}}
}
\multiput(12,16)(1,1){3}{\line(1,0){1}}      %\mu^s
\multiput(13,16)(1,1){2}{\line(0,1){1}}
\put(15,18){\line(0,1){2}}
\put(12,20){\line(1,0){4}}
\multiput(8.5,11.5)(0.5,0.75){6}{
	\makebox(0,0){$\cdot$}
}
\end{picture}
}
\end{center}
%************************** mu *********************************:
We shall consider the limit $L_J-L_{J+1} \rightarrow \infty$
($1 \le \forall J \le s,  L_{s+1}=0$) with $l_J, \mu^J$
and $L_J \equiv r_J$ mod $n$ fixed.
Our Theorem \ref{th:main_th} together with
Proposition \ref{prop:P-crystal} assert that
such limit of (\ref{eq:tensor_prod}) is bijective to
$$
\bigotimes_{J=1}^s \Bigl(
\bigoplus_{p \in {\cal H}(l_J\Lambda_{r_J}, \mu^{J})}
B\bigl(l_J\Lambda_{r_J} + af(\wt p) -
(E(p) - \overline{E}(l_J\Lambda_{r_J}, \mu^{J}))\delta \bigr) \Bigr)
$$
as affine weighted crystals.
Here $E(p)$ is the energy described in Section \ref{sect1.2},
${\cal H}(l\La_r,\mu^*)$ is
a set of restricted paths (\ref{eq:defH}) and
$\overline{E}(l\Lambda_r, \mu^*)$ is a ground state energy
(\ref{gsenergy}).
The theorem implies a ``factorization'' into a tensor product of the
pieces $J=1, \ldots, s$, and
each piece itself is a direct sum of crystals of certain integrable
highest weight modules.
In the module language it corresponds to
$$
{\cal V} = \bigotimes_{J=1}^s \Bigl(
\bigoplus_{p \in {\cal H}(l_J\Lambda_{r_J}, \mu^{J})}
V\bigl(l_J\Lambda_{r_J} + af(\wt p) -
(E(p) - \overline{E}(l_J\Lambda_{r_J}, \mu^{J}))\delta \bigr) \Bigr).
$$
See also (\ref{vcal2}).
With various choices of $\{l_J, r_J, \mu^J\}_{J=1}^s$,
this ${\cal V}$ covers a large
family of (generally reducible) $U_q(\widehat{\goth{sl}}_{\, n})$-modules.

In Section \ref{sect3} we introduce three kinds of
paths and the associated $q$-polynomials by extending
those in \cite{KMOTU3} naturally to the
inhomogeneous case.
We call them the unrestricted, classically restricted and
(level $l$) restricted one dimensional sums (1dsums) and denote by
$g_\mu(\la), X_\mu(\la)$ and $X^{(l)}_\mu(\la)$, respectively.
(Their analogues $g'_\mu(\la), X'_\mu(\la)$ and
$X^{(l)'}_\mu(\la)$ for the antisymmetric tensor case
$B_{(1^{\mu_1})}\otimes \cdots \otimes B_{(1^{\mu_m})}$
are also introduced.)
By definition they all have an expression
$\sum_p q^{E(p)}$ where the sum runs over the
weight $\la$ subset of the corresponding set of paths.
For $\mu$ finite,  Proposition \ref{pr:1dsumbyk} relates the 1dsums to
the Kostka-Foulkes polynomials as
\begin{eqnarray*}
g_\mu(\lambda) & = &\sum_{\eta (l(\eta) \le n)}K_{\eta \lambda}(1)
K_{\eta \mu}(q), \nonumber \\
{X}_\mu(\lambda) & = & K_{\lambda \mu}(q), \nonumber
\end{eqnarray*}
where the latter is due to \cite{NY}.
On the other hand, in the large $\mu$ limit
we have  Proposition \ref{1dsumlimit} as
a corollary of Theorem \ref{th:main_th}.
It identifies the limits of the 1dsums
$g_\mu(\la), X_\mu(\la)$ and $X^{(l)}_\mu(\la)$
with the string function $c^{\cal V}_\la(q)$,
the classical branching function $b^{\cal V}_\la(q)$ and
the coset branching function $a^{{\cal V}\otimes V(l_0\La_0)}_\la(q)$,
which are defined in (\ref{defc})--(\ref{defa}) and detailed in
(\ref{glimit2})--(\ref{xlimit2}).
($l_0 = l - \sum_{J=1}^s l_J$.)
Thus it is an important clue to investigate the limiting behaviour of
the Kostka-Foulkes and related polynomials
for the study of these characters.

Until the end of Section \ref{sect3} the paper only concerns
the bijection of crystals
and its general consequences on the 1dsums, which are independent of
the concrete expressions.
The rest of the paper is devoted to our second theme,
explicit formulae of the 1dsums and their limits in
{\em fermionic forms}.
By fermionic forms we roughly mean those polynomials or series which are
free of signs, admit a quasi-particle interpretation or
have an origin in the Bethe ansatz, etc.
Thanks to the absence of  signs they
are suitable for studying the limiting behaviour
and serve as a key to establish various formulae for the
characters related to the affine Lie algebras
and Virasoro algebra.
The main aim of Sections \ref{sect4} and \ref{discussion}
is to illustrate this thesis on several examples.

In Section \ref{ssec4.1} we consider the 1dsums.
As a prototype example we quote a Bethe ansatz type
fermionic formula for the Kostka-Foulkes polynomials
obtained by Kirillov and Reshetikhin \cite{KR}
in Proposition \ref{pro:ffk}.
One of our main
results in Section \ref{ssec4.1}
is the fermionic formulae for the unrestricted
1dsums $g_\mu(\la)$ and $g'_\mu(\la)$ in Propositions \ref{pro:ffkk}
and \ref{pro:ffkkp}.  The former reads
\begin{eqnarray}
g_\mu(\la)=\sum_{\eta (l(\eta) \le n)}K_{\eta\la}(1) K_{\eta\mu}(q) & = &
\sum_{\{\nu \}} q^{\phi(\{\nu \})}
\prod_{{\scriptstyle 1 \le a \le n-1} \atop
   {\scriptstyle 1 \le i \le \mu_1}}
\left[ \begin{array}{c} \nu^{(a+1)}_i -  \nu^{(a)}_{i+1}
 \\   \nu^{(a)}_i -  \nu^{(a)}_{i+1}
\end{array} \right],\label{eq:intro_ferm}\\
\phi(\{\nu \}) & = & \sum_{a=0}^{n-1} \sum_{i=1}^{\mu_1}
\left( \begin{array}{c} \nu^{(a+1)}_i -  \nu^{(a)}_i \\
2 \end{array} \right),\nonumber
\end{eqnarray}
where the sum (\ref{eq:intro_ferm}) runs over all sequences of diagrams
$\nu^{(1)},\cd,\nu^{(n-1)}$ such that
\begin{eqnarray*}
&&\emptyset =: \nu^{(0)} \subset \nu^{(1)} \subset \cdots \subset
\nu^{(n-1)} \subset \nu^{(n)}:= \mu', \nonumber \\
&&\vert \nu^{(a)} \vert = \la_1 + \cdots + \la_a \quad
\mbox{ for } 1 \le a \le n-1.\nonumber
\end{eqnarray*}
The formula (\ref{eq:intro_ferm}) is a far generalization of the
corresponding results obtained in
\cite{DJKMO1,DJKMO2,JMO,Ki,SW}.

In Section \ref{sec:fflimit} we calculate the limits
of the fermionic forms.
Proposition \ref{pro:kklim2} provides
a fermionic formula of the string function
$c^{\cal V}_\la(q)$ for the tensor product module
${\cal V} = \otimes_{J=1}^s V(l_J \Lambda_{r_J})$.
This is obtained by computing the limit of
(\ref{eq:intro_ferm}) with $\forall \mu^J = \emptyset$.
When $s=1, r_1=0$ it reduces to the one conjectured
in \cite{KNS}, announced in \cite{FS} and proved in \cite{Ge2}.
Part of some other results have also been obtained by
G. Georgiev \cite{Ge1,Ge2} and for $\widehat{\goth{sl}}_{\, 2}$ by
Schilling and Warnaar \cite{SW}.
Another important result is Proposition \ref{pro:klim},
which shows that
the limit $L \rightarrow \infty$ ($L \equiv 0$ mod $n$) of
$K_{\la (l^L)}(q)$ is expressed as a sum involving
a bilinear product of the Kostka-Foulkes polynomial and
its restricted analogue.
Under the conjecture (\ref{con:ffrk2})
this proves the
$\La = l\La_0$ case of the  spinon character formula
conjectured in \cite{NY2}:
$$
b^{V(\La)}_\la(q)
= \sum_\eta
\frac{X'_{\eta}(\la)\,
X^{(l)'}_\eta(\La)}
{(q)_{\zeta_1}\cdots (q)_{\zeta_{n-1}} },
$$
where the sum $\sum_\eta$ runs over the partitions
$\eta = \left( (n-1)^{\zeta_{n-1}},\ldots ,1^{\zeta_{1}} \right)$
satisfying $\vert \eta \vert \equiv \vert \la \vert$ mod $n$.
The numerators are the 1dsums
associated with the antisymmetric tensors defined in Section \ref{sect3}.
Proposition \ref{pro:rklim} is a similar
result related to an RSOS version of the spinon
character formula.

Section \ref{discussion} contains further generalizations.
In particular we have Conjecture \ref{con:multikns} on the
fermionic formula of the string function of the
module $V(l_1\La_0) \otimes \cdots \otimes V(l_s\La_0)$
for arbitrary non-twisted affine Lie algebra $X^{(1)}_n$.

Let us close with a few more comments
on the limiting behaviour
of the Kostka-Foulkes polynomials.
Its study was initiated by R. Gupta \cite{G}
and R. Stanley \cite{S}
(see also \cite{Ki2}) in connection to investigation of the stable
behaviour of some characters of the special linear group $SL(n)$.
In the context of integrable systems, it was initiated
by A.N. Kirillov \cite{Ki}
and continued by Nakayashiki and Yamada \cite{NY}.
The $s = 1$ case of the large $\mu$ limit
in this paper corresponds to the
so-called thermodynamical Bethe Ansatz limit \cite{Ki}, when for all
$i$ $\la_i\to\infty$ and $\mu'_i\to\infty$, but all differences
$\la_i-\la_{i+1}$ and $\mu'_i-\mu'_{i+1}$ are fixed.
Mathematically, there are yet other interesting limits.
For example, it is known \cite{KMOTU} that the 1dsums
$g'_\mu(\la), X'_\mu(\la)$ and $X^{(l)'}_\mu(\la)$
(\ref{3.4}) based on the antisymmetric tensors yield
just level 1 characters in the large
$\mu$ limit under the replacement $q \rightarrow q^{-1}$.
On the other hand from Proposition \ref{pr:1dsumbyk}
they should emerge also from the limit of
$\sum_{\eta (\eta_1 \le n)}K_{\eta' \lambda}K_{\eta \mu}(q^{-1})$
or $K_{\lambda' \mu}(q^{-1})$, etc.
Starting from their fermionic forms one can verify this
easily by using the Durfee rectangle identity at most.
Another interesting limit is the
so-called thermodynamical limit, when $\la_1\to\infty$ and $\mu'_1\to
\infty$, but $(\la_2,\cd,\la_n)$ and $(\mu'_2,\cd,\mu'_m)$ are fixed.
In this limit the Kostka-Foulkes polynomial $K_{\la\mu}(q)$ tends to
some rational function, see, e.g., \cite{Ki}.
%

%\end{document}

\section{Crystals}\label{sect1}

\subsection{Preliminaries} \label{subsec:preli}
We recapitulate necessary facts and notations concerning crystals
of the quantum affine algebra $U_q(\slchap)$. Let $\alpha_i,h_i,
\La_i$ ($i=0,1,\cd,n-1$) be the simple roots, simple coroots,
fundamental weights for the affine Lie algebra $\slchap$. For our
convenience we set $\La_{i'}=\La_i$ for any $i'\in\Z$ such that
$i'\equiv i$ mod $n$. Let $(\cdot|\cdot)$ be the standard bilinear
form normalized by $(\alpha_i|\alpha_i)=2$. The following value
will be used later: $(\La_i|\La_j)=\min(i,j)-ij/n$ ($0\le i,j<n$).
Let $\delta=\sum_{i=0}^{n-1}\alpha_i$ denote the null root, and
$c=\sum_{i=0}^{n-1}h_i$ the canonical central element. Let
$P=\oplus_{i=0}^{n-1}\Z\La_i\oplus\Z\delta$ be the weight lattice.
We define the following subsets of $P$: $P^+=\sum_{i=0}^{n-1}
\Zn\La_i$, $P^+_l=\{\la\in P^+\mid \langle\la,c\rangle=l\}$,
$\ol{P}=\sum_{i=1}^{n-1}\Z\ol{\La}_i$, $\ol{P}^+=\sum_{i=1}^{n-1}
\Zn\ol{\La}_i$. Here $\ol{\La}_i=\La_i-\La_0$ is the classical
part of $\La_i$. This map $\ol{\phantom{\La}}$ is extended to
a map on $P$ so that it is $\Z$-linear. To consider finite dimensional
$U'_q(\slchap)$-modules, the classical weight lattice $\Pcl=P/\Z\delta$
is also needed. We further define the following subsets of $\Pcl$:
$\Pcl^+=\{\la\in \Pcl\mid \langle\la,h_i\rangle\ge0\mbox{ for any }i\}$,
$\Pcll=\{\la\in \Pcl^+\mid \langle\la,c\rangle=l\}$. We introduce
an element $\La^{cl}_i\in \Pcl$ by $\La^{cl}_i=\La_i$ mod $\Z\delta$,
and fix the map $af:\Pcl\rightarrow P$ by
$af(\La^{cl}_i)=\La_i$. See Section 3.1 of \cite{KMN1} for the details
of $\Pcl,af$, etc.

The irreducible highest weight module $V(\la)$ with highest weight
$\la\in P^+$ has a crystal base $(L(\la),B(\la))$\cite{Kas}.
We denote the highest weight vector in $B(\la)$ by $u_\la$. On the crystal
$B=B(\la)$, the actions of Kashiwara operators $\et{i},\ft{i}$
($i=0,1,\cd,n-1$) are given:
\begin{eqnarray}
\et{i}&:& B\longrightarrow B\sqcup\{0\},\\
\ft{i}&:& B\longrightarrow B\sqcup\{0\}.
\end{eqnarray}
For $b,b'\in B$, $\ft{i}b=b'$ is equivalent to $b=\et{i}b'$.
Setting $\veps_i(b)=\max\{n\in\Zn\mid\et{i}^nb\neq0\}$,
$\vphi_i(b)=\max\{n\in\Zn\mid\ft{i}^nb\neq0\}$, we have
$\vphi_i(b)-\veps_i(b)=\langle h_i,\wt b\rangle$.

Crystals from the category of finite dimensional modules are also
important. Let $V_{(l)}$ be the symmetric tensor representation of
$U'_q(\slchap)$ of degree $l$. ($U'_q(\slchap)$ is the subalgebra
of $U_q(\slchap)$ generated by $e_i,f_i,q^h$ ($h\in(\Pcl)^*$).)
$V_{(l)}$ also has a crystal base $(L_{(l)},B_{(l)})$. We note that
$B_{(l)}$ is a $\Pcl$-weighted crystal. As a set, $B_{(l)}$ is
described as
\[
B_{(l)}=\{(x_1,\cd,x_n)\in\Zn^n\mid x_1+\cd+x_n=l\}.
\]
It can be identified with the set of semi-standard tableaux of shape
$(l)$ with letters from $\{1,2,\cd,n\}$. The crystal structure of
$B_{(l)}$ is given by
\begin{eqnarray}
\ft{0}(x_1,\cd,x_n)&=&(x_1+1,\cd,x_n-1),\\
\ft{i}(x_1,\cd,x_i,x_{i+1},\cd,x_n)
&=&(x_1,\cd,x_i-1,x_{i+1}+1,\cd,x_n)\quad(i\neq0).
\end{eqnarray}
If $x_i$ becomes negative upon application, $(x_1,\cd,x_n)$ should
be understood as $0$. For $b=(x_1,\cd,x_n)\in B_{(l)}$, we have
$\veps_i(b)=x_{i+1},\vphi_i(b)=x_i(i\neq0),=x_n(i=0)$ and $\wt b
=\sum_{i=1}^n x_i(\La^{cl}_i-\La^{cl}_{i-1})$.
It is sometimes convenient to write $x_i(b)$ for a component $x_i$
of $b\in B_{(l)}$.

We also review the $\Pcl$-weighted crystal $B_{(1^l)}$ of the
anti-symmetric tensor representation of $U'_q(\slchap)$ of degree $l$
($l<n$). As a set, it is described as
\[
B_{(1^l)}=\{(x_1,\cd,x_n)\in\{0,1\}^n\mid x_1+\cd+x_n=l\}.
\]
It can be identified with the set of semi-standard tableaux of
shape $(1^l)$ with letters from $\{1,2,\cd,n\}$. The crystal structure
is given similarly.

For two crystals $B_1$ and $B_2$, the tensor product $B_1\ot B_2$
is defined.
\[
B_1\ot B_2=\{b_1\ot b_2\mid b_1\in B_1,b_2\in B_2\}.
\]
The actions of $\et{i}$ and $\ft{i}$ are defined by
\begin{eqnarray}
\et{i}(b_1\ot b_2)&=&\left\{
\begin{array}{ll}
\et{i}b_1\ot b_2&\mbox{ if }\vphi_i(b_1)\ge\veps_i(b_2)\\
b_1\ot \et{i}b_2&\mbox{ if }\vphi_i(b_1) < \veps_i(b_2),
\end{array}\right.\\
\ft{i}(b_1\ot b_2)&=&\left\{
\begin{array}{ll}
\ft{i}b_1\ot b_2&\mbox{ if }\vphi_i(b_1) > \veps_i(b_2)\\
b_1\ot \ft{i}b_2&\mbox{ if }\vphi_i(b_1)\le\veps_i(b_2).
\end{array}\right.
\end{eqnarray}
Here $0\ot b$ and $b\ot0$ are understood to be $0$.
$\veps_i,\vphi_i$ and $\wt$ are given by
\begin{eqnarray}
\veps_i(b_1\ot b_2)&=&
\max(\veps_i(b_1),\veps_i(b_1)+\veps_i(b_2)-\vphi_i(b_1)),\label{eq:eps-tens}\\
\vphi_i(b_1\ot b_2)&=&
\max(\vphi_i(b_2),\vphi_i(b_1)+\vphi_i(b_2)-\veps_i(b_2)),\\
\wt(b_1\ot b_2)&=&\wt b_1+\wt b_2.
\end{eqnarray}

\subsection{Nakayashiki-Yamada's energy function}\label{sect1.2}
We review the energy function introduced by Nakayashiki and Yamada
\cite{NY}. Using this function, they represented the Kostka
polynomial $K_{\la\mu}(q)$ as a sum over paths
$b_1\ot\cd\ot b_m\in B_{(\mu_1)}\ot\cd\ot B_{(\mu_m)}$ or
$b_1\ot\cd\ot b_m\in B_{(1^{\mu_1})}\ot\cd\ot B_{(1^{\mu_m})}$
($\mu=(\mu_1,\cd,\mu_m)$) satisfying certain conditions.

Let us consider crystals $B_1$ and $B_2$ of finite dimensional
$U'_q(\slchap)$-modules. Assume the following conditions:
\begin{eqnarray}
&&B_1\ot B_2\mbox{ is connected.}\label{eq:connected}\\
&&B_1\ot B_2\mbox{ is isomorphic to }B_2\ot B_1.\label{eq:iso}
\end{eqnarray}
Suppose $b_1\ot b_2\in B_1\ot B_2$ is mapped to $b'_2\ot b'_1
\in B_2\ot B_1$ under the isomorphism. A $\Z$-valued function
$H$ on $B_1\ot B_2$ is called an {\em energy function} if for any $i$
and $b_1\ot b_2\in B_1\ot B_2$ such that $\et{i}(b_1\ot b_2)
\neq0$ it satisfies
\begin{eqnarray}
H(\et{i}(b_1\ot b_2))&=H(b_1\ot b_2)+1
&\mbox{ if }i=0,\vphi_0(b_1)\geq\veps_0(b_2),\vphi_0(b'_2)\geq\veps_0(b'_1),
\nonumber\\
&=H(b_1\ot b_2)-1
&\mbox{ if }i=0,\vphi_0(b_1)<\veps_0(b_2),\vphi_0(b'_2)<\veps_0(b'_1),
\nonumber\\
&\hspace{-6mm}=H(b_1\ot b_2)&\mbox{ otherwise}.\label{eq:e-func}
\end{eqnarray}
Explicit descriptions of the isomorphism (\ref{eq:iso}) and energy function
in the cases of $(B_1,B_2)=(B_{(k)},B_{(l)})$ and $(B_{(1^k)},B_{(1^l)})$
($k\geq l$) are given in the next subsection.

Let $B_i$ ($i=1,\cd,m$) be finite crystals such that $B_i$ and $B_j$
satisfy both (\ref{eq:connected}) and (\ref{eq:iso}) for any $i,j$ ($i<j$).
Using the isomorphism (\ref{eq:iso}), we define $b^{(i)}_j$ ($i<j$) by
\begin{eqnarray*}
&&
\begin{array}{ccccc}\hspace{-5mm}
B_i\ot\cd\ot B_{j-1}\ot B_j&\simeq&
B_i\ot\cd\ot B_j\ot B_{j-1}&\simeq&\cd\\
b_i\ot\cd\ot b_{j-1}\ot b_j&\mapsto&
b_i\ot\cd\ot b^{(j-1)}_j\ot b'_{j-1}&\mapsto&\cd
\end{array}\\
&&\hspace{5cm}
\begin{array}{ccc}
\cd&\simeq&B_j\ot B_i\ot\cd\ot B_{j-1}\\
\cd&\mapsto&b^{(i)}_j\ot b'_i\ot\cd\ot b'_{j-1},
\end{array}
\end{eqnarray*}
and set $b^{(i)}_i=b_i$. Consider an element $p=b_1\ot\cd\ot b_m$
of $B_1\ot\cd\ot B_m$.
We call the following quantity the {\em energy} of $p$.
\[
E(p)=\sum_{i<j}H_{ij}(b_i\ot b^{(i+1)}_j).
\]
Here $H_{ij}$ is the energy function on $B_i\ot B_j$ defined previously.

Consider the case when $B_i=B_{(\mu_i)}$ with $\mu=(\mu_1,\cd,\mu_m)$
a partition. The set of highest weight crystals
of weight $\la$ in $B_{(\mu_1)}\ot\cd\ot B_{(\mu_m)}$ with respect to
$U_q(\goth{sl}_{\,n})$ is known to be bijective to the set of
semi-standard tableaux of shape $\la$ and weight $\mu$. It was shown
in \cite{NY} that the energy of $p$ coincides with the charge of the
corresponding tableau in the sense of Lascoux-Sch\"utzenberger \cite{LS}.

Let $p$ be as above, and consider the following condition for $b_1\in B_1$.
\begin{eqnarray}
&&\mbox{For any $j(\neq1)$ and $b_j\in B_j$,}\nonumber\\
&&\mbox{if $\et{0}(b_1\ot b_j)=b_1\ot\et{0}b_j$,
then $\et{0}(b'_j\ot b'_1)=b'_j\ot\et{0}b'_1$.}\label{eq:cond_b1}
\end{eqnarray}
Here $b_1\ot b_j$ is mapped to $b'_j\ot b'_1$ under the isomorphism
$B_1\ot B_j\simeq B_j\ot B_1$. We have the following
representation-theoretic interpretation of energy.

\begin{proposition} \label{prop:interpre}
Let $p$ be as above. If $i\neq0$ and $\et{i}p\neq0$, then
\[
E(\et{i}p)=E(p).
\]
If $\et{0}p=b_1\ot\cd\ot\et{0}b_k\ot\cd\ot b_m\neq0$ with $k\neq1$
and $b_1$ satisfies the condition (\ref{eq:cond_b1}), then
\[
E(\et{0}p)=E(p)-1.
\]
\end{proposition}
The case of $i\neq0$ is clear. The case of $i=0$ reduces to the following.

\begin{lemma}
If $\et{0}p=b_1\ot\cd\ot\et{0}b_k\ot\cd\ot b_m\neq0$ with $k\neq1$
and $b_1$ satisfies the condition (\ref{eq:cond_b1}), then
\begin{eqnarray}
E^{(j)}(\et{0}p)&=E^{(j)}(p)-1&(j=k)\label{eq:E(j)-a}\\
&\hspace{-6mm}=E^{(j)}(p)&(j\neq k).\label{eq:E(j)-b}
\end{eqnarray}
Here $E^{(j)}(p)=\sum_{i=1}^{j-1}H_{ij}(b_i\ot b^{(i+1)}_j)$.
\end{lemma}

\Proof
Set $\btilde_i=\et{0}b_i(i=k),b_i(i\neq k)$. In the case of $j<k$,
we have $\btilde_i=b_i,\btilde^{(i+1)}_j=b^{(i+1)}_j$ for $1\leq i\leq j-1$.
This shows (\ref{eq:E(j)-b}) when $j<k$.

To show it in the case of $j\geq k$, we rewrite (\ref{eq:e-func}) in
the following manner.
\begin{eqnarray*}
H(\et{0}(b_1\ot b_2))&=H(b_1\ot b_2)+1
&\mbox{ if }\et{0}(b_1\ot b_2)=\et{0}b_1\ot b_2\mbox{ and }\\
&&\hspace{4.6mm}\et{0}(b'_2\ot b'_1)=\et{0}b'_2\ot b'_1,\\
&=H(b_1\ot b_2)-1
&\mbox{ if }\et{0}(b_1\ot b_2)=b_1\ot\et{0}b_2\mbox{ and }\\
&&\hspace{4.6mm}\et{0}(b'_2\ot b'_1)=b'_2\ot\et{0}b'_1,\\
&\hspace{-6mm}=H(b_1\ot b_2)&\mbox{ otherwise}.
\end{eqnarray*}
Let $k'$ be the largest integer such that
\begin{eqnarray*}
\begin{array}{c}
B_1\ot\cd\ot B_k\ot\cd\ot B_j\\
b_1\ot\cd\ot\et{0}b_k\ot\cd\ot b_j
\end{array}&&\\
&&\\
&&\hspace{-3cm}
\begin{array}{cc}
\simeq&
B_1\ot\cd\ot B_k\ot B_j\ot\cd\ot B_{j-1}\\
\mapsto&
b_1\ot\cd\ot\et{0}b_k\ot b^{(k+1)}_j\ot\cd\ot b'_{j-1}\\
&\\
\simeq&
B_1\ot\cd\ot B_j\ot B_k\ot\cd\ot B_{j-1}\\
\mapsto&
b_1\ot\cd\ot\et{0}b^{(k)}_j\ot b'_k\ot\cd\ot b'_{j-1}\\
&\\
\simeq&
B_1\ot\cd\ot B_{k'}\ot B_j\ot\cd\ot B_{j-1}\\
\mapsto&
b_1\ot\cd\ot b_{k'}\ot\et{0}b^{(k'+1)}_j\ot\cd\ot b'_{j-1}\\
&\\
\simeq&
B_1\ot\cd\ot B_j\ot B_{k'}\ot\cd\ot B_{j-1}\\
\mapsto&
b_1\ot\cd\ot b^{(k')}_j\ot\et{0}b'_{k'}\ot\cd\ot b'_{j-1}.
\end{array}
\end{eqnarray*}
Note that $1\leq k'<k$ if $j=k$ and $1\leq k'\leq k$ if $j>k$.
The existence of such $k'$ is guaranteed
by (\ref{eq:cond_b1}). If $k'=k$, the 3rd and 4th terms should be omitted.
{}From the above property of $H$, we get the following.
If $j=k$, $H_{ij}(\btilde_i,\btilde^{(i+1)}_j)-H_{ij}
(b_i,b^{(i+1)}_j)=-1\,(i=k'),=0\,(\mbox{otherwise})$.
If $j>k$ and $k'\neq k$, $H_{ij}(\btilde_i,\btilde^{(i+1)}_j)-
H_{ij}(b_i,b^{(i+1)}_j)=1\,(i=k),-1\,(i=k'),=0\,
(\mbox{otherwise})$.
If $j>k$ and $k'=k$, $H_{ij}(\btilde_i,\btilde^{(i+1)}_j)-
H_{ij}(b_i,b^{(i+1)}_j)=0$. These imply (\ref{eq:E(j)-a}) and
(\ref{eq:E(j)-b}) respectively.
\qed

\subsection{Evaluation of the energy} \label{subsec:energy}
Here we give an explicit procedure to obtain
the energy function $H : B_1 \ot B_2 \rightarrow {\bf Z}$
and the isomorphism $\iota : B_1 \ot B_2 \rightarrow B_2 \ot B_1$
in the case of symmetric tensor representations
$(B_1,B_2)=(B_{(k)},B_{(l)})$ ($k \geq l$).

Let $b_1 \ot b_2$ be an element in $B_1 \ot B_2$
such as $b_1=(x_1,\ldots,x_n)$ and $b_2=(y_1,\ldots,y_n)$.
We represent $b_1 \ot b_2$ by the two column diagram.
Each column has $n$ rows, enumerated as 1 to $n$ from the top to the bottom.
We put $x_i$ (resp. $y_i$) dots $\bullet$ in the $i$-th row
of the left (resp. right) column.
The labels of rows are sometimes extended to ${\bf Z}/n {\bf Z}$.

\begin{proposition}\label{prop:H-rule}
The rule to obtain the energy function $H$ and the
isomorphism $\iota$ is as follows.
\begin{itemize}
\item[(1)]
Pick any dot, say $\bullet_a$, in the right column and connect it
with a dot $\bullet_a'$ in the left column by a line
(which we call $H$-line).
The partner $\bullet_a'$ is chosen
{}from the dots which are in the lowest row among all dots
whose positions are higher than that of $\bullet_a$.
If there is no such dot, we return to the bottom and
the partner $\bullet_a'$ is chosen from the dots
in the lowest row among all dots.
In the latter case, we call such a pair or line ``winding".
% ********************* Figure (H-line Proposition) ***********
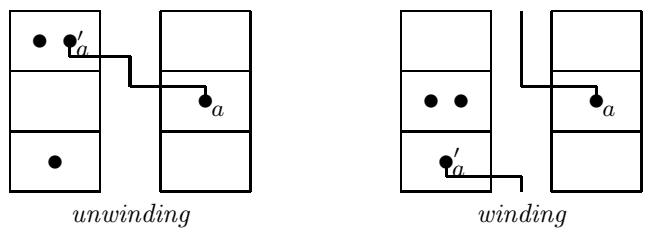
\begin{center}
\unitlength=0.8mm
\begin{picture}(105,35)(0,-5)
\multiput(0,0)(65,0){2}{
	\multiput(0,0)(25,0){2}{
		\multiput(0,0)(0,10){4}{\line(1,0){15}}
		\multiput(0,0)(15,0){2}{\line(0,1){30}}
	}
	\put(32.5,15){\circle*{2}}
	\put(33.5,14.3){\makebox(0,0)[lt]{\small $a$}}
}
\put( 5,25){\circle*{2}}
\put(10,25){\circle*{2}}
\put(11.2,26.8){\makebox(0,0)[lt]{$'$}}
\put(11,24.3){\makebox(0,0)[lt]{\small $a$}}
\put(7.5,5){\circle*{2}}
\put(70,15){\circle*{2}}
\put(75,15){\circle*{2}}
\put(72.5,5){\circle*{2}}
\put(73.7,7.3){\makebox(0,0)[lt]{$'$}}
\put(73.5,4.3){\makebox(0,0)[lt]{\small $a$}}
\thicklines                        % H-lines
\put(20,17.5){\line(1,0){12.5}}
\put(10,22.5){\line(1,0){10}}
\put(20,22.5){\line(0,-1){5}}
\put(10,25){\line(0,-1){2.5}}
\put(32.5,15){\line(0,1){2.5}}
\put(  85,17.5){\line(1,0){12.5}}
\put(72.5, 2.5){\line(1,0){12.5}}
\put(  85, 2.5){\line(0,-1){2.5}}
\put(72.5,   5){\line(0,-1){2.5}}
\put(  85,17.5){\line(0,1){12.5}}
\put(97.5,  15){\line(0,1){2.5}}
\put(20,-2){\makebox(0,0)[t]{unwinding}}
\put(85,-2){\makebox(0,0)[t]{winding}}
\end{picture}
\end{center}
% ********************* Figure (H-line Proposition) ***********
\item[(2)]
Repeat the procedure (1) for the remaining unconnected dots
$(l-1)$-times.
\item[(3)]
The isomorphism $\iota$ is obtained by sliding
the remaining $(k-l)$ unpaired dots in the left column to
the right.
\item[(4)]
The value of the energy function is the number of the ``winding" pairs.
\end{itemize}
\end{proposition}

The $H$ and $\iota$ obtained by this rule have the correct property
as the energy function and isomorphism.
This fact has been proved in \cite{NY} Section 3,
where the following two lemmas are also proved.

\begin{lemma}
The map $\iota$ determined by the above rule
is independent of the order of drawing lines.
\end{lemma}

\Proof
Suppose there exists a dot (say $\bullet_{L1}$ ) in
the left column that is unpaired in one ordering (say A)
and paired in another (say B).
Let $\bullet_{R1}$ be the partner of $\bullet_{L1}$ in B, and
let $\bullet_{L2}$ be the partner of $\bullet_{R1}$ in A.
Then $\bullet_{L2}$ must be paired with some dot (say $\bullet_{R2}$ )
in B since the $R_1L_1$ line already passes through $\bullet_{L2}$.
This process of determining $R_1,L_2,R_2,\cdots$ does not stop and
$R_1,R_2,\cdots$ are all distinct.
This is a contradiction, since the number of dots in the right
column is finite (=$l$).
Therefore the set of end points of the $H$-lines is independent of the
order in which they are drawn and the lemma is proved. \qed

\begin{lemma}
The value of the function $H$ determined by the above rule
is independent of the order of drawing lines.
\end{lemma}

\Proof
For $j$ ($1 \leq j \leq n$) we assign a non-negative
integer $h_j(A)$ as the number of lines passing the $j$-th row.
Here the word ``passing" is defined as follows.
Let $\alpha$ be a line starting from the $i$-th row and
ending at the $j$-th row, then the line $\alpha$ passes
the $k$-th row if and only if
(1) $i>k>j$ (for non-winding $\alpha$) or
(2) $k<i$ or $k>j$ (for winding $\alpha$).
$h_j(A)$ may depend on the order $A$ of drawing $H$-lines.
$h_j(A)$ is subject to the following relation:
\begin{equation}
h_{j+1}(A)=h_{j}(A)+e_j-s_{j+1},
\label{eq:h-number}
\end{equation}
where $e_j$ is the number of dots in the left column
which are end points of the $H$-lines and sitting in the $j$-th row.
Similarly, $s_j$ is the number of dots in the right column
(starting points of the $H$-lines) sitting in the $j$-th row.
% ********************* Figure (H-line Lemma) ***********
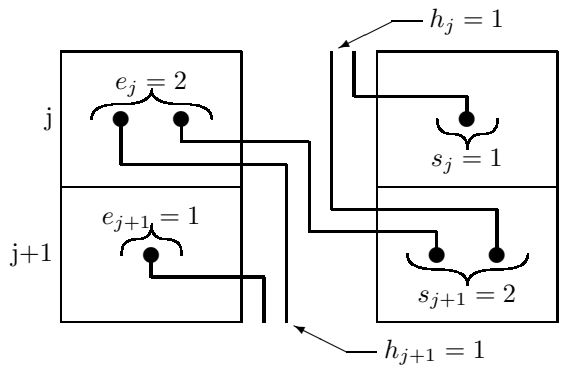
\begin{center}
\unitlength=1mm
\begin{picture}(73,43)(-7,0)
\multiput(0,0)(42,0){2}{
	\multiput(0,0)(0,18){3}{\line(1,0){24}}
	\multiput(0,0)(24,0){2}{\line(0,1){36}}
}
\put(12,9){\circle*{2}}
\multiput(50,9)(8,0){2}{\circle*{2}}
\multiput(8,27)(8,0){2}{\circle*{2}}
\put(54,27){\circle*{2}}
\thicklines                        % H-lines
\put(39,36){\line(0,-1){6}}
\put(54,30){\line(-1,0){15}}
\put(54,30){\line(0,-1){3}}
\put(36,36){\line(0,-1){21}}
\put(58,15){\line(-1,0){22}}
\put(58,15){\line(0,-1){6}}
\put(16,27){\line(0,-1){3}}
\put(33,24){\line(-1,0){17}}
\put(33,24){\line(0,-1){12}}
\put(50,12){\line(-1,0){17}}
\put(50,12){\line(0,-1){3}}
\put( 8,27){\line(0,-1){6}}
\put(30,21){\line(-1,0){22}}
\put(30,21){\line(0,-1){21}}
\put(12, 9){\line(0,-1){3}}
\put(27, 6){\line(-1,0){15}}
\put(27, 6){\line(0,-1){6}}
\thinlines
\put(-1,9){\makebox(0,0)[r]{j+1}}
\put(-1,27){\makebox(0,0)[r]{j}}
\bezier{200}(46,9)(46,7)(50,7)
\bezier{200}(50,7)(54,7)(54,5)
\bezier{200}(58,7)(54,7)(54,5)
\bezier{200}(62,9)(62,7)(58,7)
\put(54,5){\makebox(0,0)[t]{$s_{j+1}=2$}}
\bezier{200}(50,27)(50,25)(52,25)
\bezier{200}(52,25)(54,25)(54,23)
\bezier{200}(56,25)(54,25)(54,23)
\bezier{200}(58,27)(58,25)(56,25)
\put(54,23){\makebox(0,0)[t]{$s_{j}=1$}}
\bezier{200}( 8, 9)( 8,11)(10,11)
\bezier{200}(10,11)(12,11)(12,13)
\bezier{200}(14,11)(12,11)(12,13)
\bezier{200}(16, 9)(16,11)(14,11)
\put(12,13){\makebox(0,0)[b]{$e_{j+1}=1$}}
\bezier{200}( 4,27)( 4,29)( 8,29)
\bezier{200}( 8,29)(12,29)(12,31)
\bezier{200}(16,29)(12,29)(12,31)
\bezier{200}(20,27)(20,29)(16,29)
\put(12,31){\makebox(0,0)[b]{$e_{j}=2$}}
\put(44,40){\vector(-2,-1){7}}
\put(44,40){\line(1,0){4}}
\put(49,40){\makebox(0,0)[l]{$h_j=1$}}
\put(38,-4){\vector(-2,1){7}}
\put(38,-4){\line(1,0){4}}
\put(43,-4){\makebox(0,0)[l]{$h_{j+1}=1$}}
\end{picture}
\end{center}
\vspace{0.5cm}
% ********************* Figure (H-line Lemma) ***********
The value of the energy function $H$ with respect to the order $A$
is given by $h_1(A)+s_1$.
Note that the set of end points of the $H$-lines is independent of the
order $A$, hence $e_j$ and $s_j$ are independent of the order $A$.

We prove that $h_j(A)$ does not depend on $A$.
To this end, it is sufficient to prove that for any $A$ there exists
$j$ such that $h_j(A)=0$.
In fact, for two orders $A$ and $B$,
there exists an integer $m$ such that
$$
h_j(A)=h_j(B)+m
$$
for any $j$ by (\ref{eq:h-number}).
By exchanging $A$ and $B$ if necessary, we can assume $m \geq 0$.
The existence of $j$ for which $h_j(A)=0$ means $m=0$,
since $h_j(B) \geq 0$.
The existence of such $j$ can be proved as follows.
Suppose that such $j$ does not exist.
Then there exists a sequence of $H$-lines
$\alpha_1,\alpha_2,\ldots,\alpha_k(=\alpha_0)$ such that
the end point of $\alpha_i$ is passed through by the
line $\alpha_{i-1}$.
If such a situation occurs, however, the order of drawing $H$-lines should
satisfy
$$
\hbox{order of $\alpha_1$} >
\hbox{order of $\alpha_2$} >
\ldots >
\hbox{order of $\alpha_k$} >
\hbox{order of $\alpha_1$}
$$
which is a contradiction. \qed

The following is just a corollary of Proposition \ref{prop:H-rule}.
\begin{proposition}\label{prop:E-line}
Let $\mu=(\mu_1,\ldots,\mu_m)$ be a partition and let
$p=b_1 \ot \ldots \ot b_m$ be a path in
$B_{(\mu_1)} \ot \ldots \ot B_{(\mu_m)}$.
The rule to evaluate the energy
\begin{eqnarray*}
E(p)&=&\sum_{j=1}^m E^{(j)}(p),\\
E^{(j)}(p)&=&\sum_{i=1}^{j-1}H_{(\mu_i)(\mu_j)}(b_i\ot b^{(i+1)}_j)
\quad(1\le j\le m)
\end{eqnarray*}
is given as follows.
\begin{itemize}
\item[(1)]
Pick any dot in $b_m$. According to the $H$-line rule,
connect the dot in $b_m$ with a dot in $b_{m-1}$ and
connect the dot in $b_{m-1}$ with a dot in $b_{m-2}$ and
continue this until we come to a dot in $b_1$.
We call this line $E$-line.
\item[(2)]
Repeat the procedure (1) for the remaining unconnected dots
$(\mu_m-1)$-times.
\item[(3)]
Forget about the connected dots and repeat the procedures (1) and (2)
{}from a rightmost unconnected dot. Eventually, all the dots in $p$
are decomposed into a disjoint union of $E$-lines. We call it $E$-line
decomposition of $p$.
\item[(4)]
$E^{(j)}(p)$ is given as the sum of winding numbers between $b_1$
and $b_j$ of all the $E$-lines starting from $b_i$ with $i\ge j$.
\end{itemize}
\end{proposition}

\begin{example}
Let $n=3$ and $p=(1,0,2)\otimes(0,2,0)\otimes(0,1,1)\otimes(0,1,0)
\in B_{(3)} \otimes B_{(2)} \otimes B_{(2)} \otimes B_{(1)}$.
The energy $E(p)=0+1+2+1=4$ of this path $p$ is evaluated by the following
diagram.

\begin{center}
\setlength{\unitlength}{0.7mm}
\begin{picture}(80,40)(0,-3)
\put(0,0){\line(1,0){80}}
\put(0,10){\line(1,0){80}}
\put(0,20){\line(1,0){80}}
\put(0,30){\line(1,0){80}}
\put(0,0){\line(0,1){30}}
\put(20,0){\line(0,1){30}}
\put(40,0){\line(0,1){30}}
\put(60,0){\line(0,1){30}}
\put(80,0){\line(0,1){30}}
\put(7,5){\circle*{2}}
\put(13,5){\circle*{2}}
\put(10,25){\circle*{2}}
\put(27,15){\circle*{2}}
\put(33,15){\circle*{2}}
\put(50,5){\circle*{2}}
\put(50,15){\circle*{2}}
\put(70,15){\circle*{2}}
\put(70,15){\line(0,1){10}}
\put(70,25){\line(-1,0){20}}
\put(50,25){\line(0,1){10}}
\put(50,-5){\line(0,1){10}}
\put(50,5){\line(-1,0){17}}
\put(33,5){\line(0,1){10}}
\put(33,15){\line(0,1){3}}
\put(33,18){\line(-1,0){23}}
\put(10,18){\line(0,1){7}}
\put(50,15){\line(0,1){7}}
\put(50,22){\line(-1,0){23}}
\put(27,22){\line(0,1){13}}
\put(27,-5){\line(0,1){20}}
\put(27,15){\line(-1,0){20}}
\put(7,15){\line(0,1){20}}
\put(7,-5){\line(0,1){10}}
\end{picture}
\end{center}
\end{example}

\begin{remark}
Let $p$ be as in Proposition \ref{prop:E-line}. We can check
any $b_1\in B_{(\mu_1)}$ satisfies the condition (\ref{eq:cond_b1}).
Just note that $\et{0}(b_1\ot b_j)=b_1\ot\et{0}b_j$ implies
$x_n(b_1)<x_1(b_j)$, which is sufficient for $x_n(b'_j)<x_1(b'_1)$
by Proposition \ref{prop:H-rule}.
\end{remark}
\begin{remark}
For the anti-symmetric case
$(B_1,B_2)=(B_{(1^k)},B_{(1^l)})$ $(k\geq l)$,
the rule is almost the same as in the symmetric case.
The only differences are in (1) and (4).
\begin{itemize}
\item[(1)]
The partner $\bullet_a'$ is a dot that
has the highest position among all dots whose positions are
not higher than $\bullet_a$.
If there is no such dot, we return to the top.
We call such a pair ``winding".
\item[(4)]
The value of the energy function is
$(-1)$ times the number of ``winding" pairs.
\end{itemize}
The independence of this rule with respect to the order
of drawing $H$-lines can be proved similarly.
There also exists a similar rule to obtain the energy $E(p)$
for the anti-symmetric case.
\end{remark}
\begin{remark}
The energy functions $H_{(\mu_i) (\mu_j)}$ and
$H_{(1^{\mu_i}) (1^{\mu_j})}$ can also be evaluated
by using the ``nonmovable tableaux'' as in \cite{KKN}.
\end{remark}

\section{Tensor product of crystals}\label{sect2}
In this section, we consider $B(\la)$ as a $\Pcl$-weighted crystal
except Section \ref{subsec:aff-wt}.

\subsection{Decomposition of $B(\la)\ot B_{(l)}$}
The crystal $B_{(l)}$ is known to be {\em perfect} of level $l$
\cite{KMN1,KMN2}. It means that for any $\la\in\Pcll$, we have an isomorphism
\[
B(\la)\ot B_{(l)}\simeq B(\la')
\]
for some $\la'\in\Pcll$. If the level $k$ of $\la$ is greater than $l$,
it is known that $B(\la)\ot B_{(l)}$ decomposes into a disjoint union
of crystals $B(\mu)$ with $\mu$ being a dominant integral weight of
level $k$. (See the beginning of Section 5 of \cite{KK}.) More precisely,
we have

\begin{theorem}
Let $\la\in(\Pcl^+)_k$. If $k\ge l$, then
\[
B(\la)\ot B_{(l)}\simeq\bigoplus_{b\in B^{\le\la}_{(l)}}B(\la+\wt b),
\]
where
\[
B^{\le\la}_{(l)}=\{b\in B_{(l)}\mid\veps_i(b)\le\langle h_i,\la\rangle
\mbox{ for all }i\}.
\]
\end{theorem}

\begin{corollary}
Let $\mu=(\mu_1,\cd, \mu_m)$ be a partition and $\la\in\Pcll$. If
$l\ge\mu_1$, then
\begin{equation}
B(\la)\ot B_{(\mu_1)}\ot\cd\ot B_{(\mu_m)}\simeq
\bigoplus_{p\in{\cal H}(\la,\mu)}B(\la+\wt p),
\label{eq:decomp}
\end{equation}
where
\begin{equation}
{\cal H}(\la,\mu)=\left\{p=b_1\ot\cd\ot b_m\left|
\begin{array}{l}
b_j\in B_{(\mu_j)}\mbox{ for }j=1,\cd,m,\\
\veps_i(b_j)\le\langle h_i,\la+\wt b_1+\cd+\wt b_{j-1}\rangle\\
\hspace{1.8cm}\mbox{ for all }i\mbox{ and }j=1,\cd,m
\end{array}
\right.\right\}.
\label{eq:defH}
\end{equation}
\end{corollary}
Note that $B^{\le\la}_{(l)}={\cal H}(\la,(l))$.
We define $\B(\la,\mu)$ by the RHS of (\ref{eq:decomp}).

\begin{example} \label{ex:calB}
We give examples of $\B(\la,\mu)$.
\begin{itemize}
\item[(1)]
Let $\la$ be of level $l$ and $\mu=(l^m)$. Then the perfectness of
$B_{(l)}$ fixes $p$ to be unique, and we have
\[
\B(\la,(l^m))=B(\sigma^m(\la)).
\]
Here $\sigma$ is the automorphism on $\Pcl^+$ defined by
\[
\sigma(a_0\La^{cl}_0+a_1\La^{cl}_1+\cd+a_{n-1}\La^{cl}_{n-1})
=a_{n-1}\La^{cl}_0+a_0\La^{cl}_1+\cd+a_{n-2}\La^{cl}_{n-1}.
\]
\item[(2)]
Consider the case when $\la=l\La^{cl}_r,\mu=(s)\,(0\le r\le n-1,l\ge s)$.
In this case, we have
\begin{eqnarray*}
{\cal H}(l\La^{cl}_r,(s))&=&\{(0,\cd,0,\mathop{s}^{r+1},0,\cd,0)\},\\
\B(l\La^{cl}_r,(s))&=&B((l-s)\La^{cl}_r+s\La^{cl}_{r+1}).
\end{eqnarray*}
\item[(3)]
Consider the case when $\la=l\La^{cl}_0,\mu=(s,t)\,(l\ge s\ge t)$.
In this case, we have
\begin{eqnarray*}
{\cal H}(l\La^{cl}_0,(s,t))&=&\left\{(s,0,\cd,0)\ot b\left|
\begin{array}{l}
x_1(b)\le l-s,x_2(b)\le s,\\
x_i(b)=0\,(i\ge3)
\end{array}
\right.\right\},\\
\B(l\La^{cl}_0,(s,t))&=&\bigoplus_{i=0}^r
B((l-s-i)\La^{cl}_0+(s-t+2i)\La^{cl}_1+(t-i)\La^{cl}_2),
\end{eqnarray*}
where $r=\min(l-s,t)$.
\end{itemize}
\end{example}

We give a
characterization of the set ${\cal H}(l\La^{cl}_r,\mu)\,(0\le r<n)$,
which will often appear in later sections.
Each path $p=b_1\ot\cd\ot b_m$ in ${\cal H}(l \La^{cl}_r,\mu)$
corresponds to a sequence of Young diagrams $\nu^{(a)}$ ($a=0,\cd,m$)
such that
\begin{itemize}
\item[(1)] $\nu^{(0)}=(l^r)$,

\item[(2)] $\nu^{(a)}/\nu^{(a-1)}$ is horizontal strip of length $\mu_a$ for
           ($a=1,\cd,m$),

\item[(3)] depth of $\nu^{(a)} \leq n$ and $\nu^{(a)}_1-\nu^{(a)}_n \leq l$.

\end{itemize}
The correspondence is given by $x_i(b_a)=$ the number of nodes in the
$i$-th row of $\nu^{(a)}/\nu^{(a-1)}$.
Under this correspondence, we have $l\La^{cl}_r+\wt p=l\La^{cl}_0
+\sum_{k\ge1}{\nu^{(m)}}'_k(\La^{cl}_k-\La^{cl}_0)$, where ${\nu^{(m)}}'$
denotes the transpose of $\nu^{(m)}$.

\subsection{Limit} \label{subsec:limit}
For given partitions $\la$ and $\mu$, we define two operations
$\la\cup\mu$ and $\la+\mu$. Let $\la=(\la_1,\cd,\la_l)$,
$\mu=(\mu_1,\cd,\mu_m)$. Assuming $\la_l\ge\mu_1$, $\la\cup\mu$
is defined by $(\la_1,\cd,\la_l,\mu_1,\cd,\mu_m)$. Appending
0's if necessary, we now assume $l=m$. Then $\la+\mu$ is
defined by $(\la_1+\mu_1,\cd,\la_l+\mu_l)$. Associated to a
partition $\mu=(\mu_1,\cd,\mu_m)$, we next define a special
element $\pbar$ of $B_{(\mu_1)}\ot\cd\ot B_{(\mu_m)}$ called the
{\em ground state path}. It is given by
\begin{eqnarray}
\pbar&=&\bbar_1\ot \bbar_2\ot\cd\ot \bbar_m,\\
\bbar_j&=&(0,\cd,0,\mathop{\mu_j}^r,0\cd,0)\in B_{(\mu_j)},
\end{eqnarray}
where the position $r$ is determined by $r\equiv j$ (mod $n$),
$1\le r\le n$.

Let $l_J,L_J$ be positive integers and $\mu^J$ be a partition
such that $l_J\ge\mu^J_1$ for $J=1,\cd,s$. It is convenient to
set $l=\sum_{J=1}^sl_J$ and $L_{s+1}=0$.
We consider the following partition.
\begin{eqnarray}
\mu&=&\mu((l_1,L_1,\mu^1),\cd,(l_s,L_s,\mu^s))\\
&=&((l_1^{L_1})\cup\mu^1)+\cd+((l_s^{L_s})\cup\mu^s). \label{eq:mu}
\end{eqnarray}
%************************** mu *********************************:
\begin{center}
\unitlength=2mm
{\small
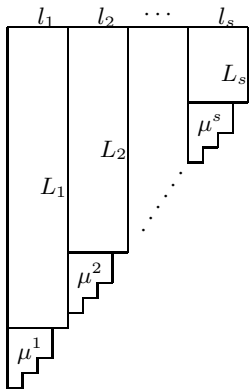
\begin{picture}(16,28)(0,1)
\put(0,25){\line(1,0){16}}
\put(0,25){\line(0,-1){24}}
\put(4,25){\line(0,-1){20}}
\put(8,25){\line(0,-1){15}}
\put(12,25){\line(0,-1){9}}
\put(16,25){\line(0,-1){5}}
\put(2,25.2){\makebox(0,0)[lb]{$l_1$}}
\put(6,25.2){\makebox(0,0)[lb]{$l_2$}}
\put(10,25.5){\makebox(0,0)[b]{$\cdots$}}
\put(14,25.2){\makebox(0,0)[lb]{$l_s$}}
\put(4,15){\makebox(0,0)[rt]{$L_1$}}
\put(8,17.5){\makebox(0,0)[rt]{$L_2$}}
\put(16,22.5){\makebox(0,0)[rt]{$L_s$}}
\put(1.5,3.5){\makebox(0,0){$\mu^1$}}
\put(5.5,8.5){\makebox(0,0){$\mu^2$}}
\put(13.5,18.5){\makebox(0,0){$\mu^s$}}
\multiput(0,1)(4,5){2}{                      %\mu^1,\mu^2
	\multiput(0,0)(1,1){3}{\line(1,0){1}}
	\multiput(1,0)(1,1){2}{\line(0,1){1}}
	\put(3,2){\line(0,1){2}}
	\put(0,4){\line(1,0){4}}
}
\multiput(12,16)(1,1){3}{\line(1,0){1}}      %\mu^s
\multiput(13,16)(1,1){2}{\line(0,1){1}}
\put(15,18){\line(0,1){2}}
\put(12,20){\line(1,0){4}}
\multiput(8.5,11.5)(0.5,0.75){6}{
	\makebox(0,0){$\cdot$}
}
\end{picture}
}
\end{center}
%************************** mu *********************************:
We would like to consider the limit when all $L_J-L_{J+1}$
($J=1,\cd,s$) go to infinity with the residue of $L_J$ mod $n$ fixed.
For $r_J$ ($0\le r_J<n,1\le J\le s$), let us define
$\P_{<\infty}((l_1,r_1,\mu^1),\cd,(l_s,r_s,\mu^s))$ to be the set
of elements $p$ such that
\begin{itemize}
\item[(1)] $p$ is an element of the limit of
           $B_{(\mu_1)}\ot\cd\ot B_{(\mu_m)}$ when all $L_J-L_{J+1}$
           ($J=1,\cd,s$) go to infinity with the residue of
           $L_J$ mod $n$ fixed as $r_J$.
\item[(2)] $E(p)-E(\pbar)$ is finite.
\end{itemize}

\begin{theorem} \label{th:main_th}
With the definitions as above, we have the following statements.
\begin{itemize}
\item[(1)] There exists a bijection
\[
\P_{<\infty}((l_1,r_1,\mu^1),\cd,(l_s,r_s,\mu^s))
\simeq\bigotimes_{J=1}^s
\P_{<\infty}((l_J,r_J,\mu^J)).
\]
\item[(2)] If $p$ is mapped to $p^1\ot\cd\ot p^s$ under this bijection,
we have
\[
E(p)-E(\pbar)=\sum_{J=1}^s(E(p^J)-E(\pbar^J)),
\]
where $\pbar^J$ denotes the ground state path of
$\P_{<\infty}((l_J,r_J,\mu^J))$.
\item[(3)] $\P_{<\infty}((l_J,r_J,\mu^J))$ is identified with the crystal
$\B(l_J\La^{cl}_{r_J},\mu^J)$.
\end{itemize}
\end{theorem}

\subsection{Proof}
Before presenting lemmas, we prepare some notations. Let $p$ be a path
corresponding to $\mu$, i.e., $p\in B_{(\mu_1)}\ot\cd\ot B_{(\mu_m)}$.
Let partitions $\mu^1,\mu^2$ be such that $\mu=\mu^1+\mu^2$. $(p^1,p^2)$
is said to be a $(\mu^1,\mu^2)$-{\em partition} of $p$ if $p^J\in
B_{(\mu^J_1)}\ot\cd\ot B_{(\mu^J_m)}$ for $J=1,2$ and $x_i(b_j)=
x_i(b^1_j)+x_i(b^2_j)$ for $i=1,\cd,n$ and $j=1,\cd,m$. Here $b_j$ and
$b^J_j$ are defined by $p=b_1\ot\cd\ot b_m$ and $p^J=b^J_1\ot\cd\ot b^J_m$.
A $(\mu^1,\cd,\mu^s)$-partition $(p^1,\cd,p^s)$ of $p$ is defined similarly.

We prepare five lemmas.

\begin{lemma} \label{lem:1}
Let $\mu^1,\mu^2,\nu$ be partitions such that $\mu^1=(l_1^{L_1})\cup\nu,
l_1\ge\nu_1$ and the depth of $\mu^2$ is smaller than $L_1$. Let $p$ be
a path corresponding to $\mu=\mu^1+\mu^2$. Then we have
\[
E(p)\ge\min(E(p^1)+E(p^2)),
\]
where the minimum is taken over all $(\mu^1,\mu^2)$-partitions
$(p^1,p^2)$ of $p$.
\end{lemma}

\Proof
First fix any $E$-line decomposition $\ell$ of $p$. Define $\ell_1$ to
be the set of $E$-lines penetrating the $L_1$-th tensor component, and
set $\ell_2=\ell\setminus\ell_1$. Next define $p^J_0=b^J_1\ot b^J_2\ot\cd$
($J=1,2$) in such a way that $x_k(b^J_j)$ is the number of dots on
$\ell_J$ at position $k$ in the $j$-th tensor component. Then we have
\[
E(p)=E(p^1_0)+E(p^2_0)\ge\min_{(p^1,p^2)}(E(p^1)+E(p^2)).
\]
Here the first equality is due to Proposition \ref{prop:E-line}.
\qed

In what follows, we set $\delta^{(n)}_{ij}=1$ ($i\equiv j$ mod $n$),
$=0$ (otherwise).

\begin{lemma} \label{lem:2}
Let $\mu^1,\mu^2,\mu$ and $p$ be as above. Let $(p^1,p^2)$ be
a $(\mu^1,\mu^2)$-partition of $p$. We further assume that by setting
$p^1=b^1_1\ot b^1_2\ot\cd$, we have $x_i(b^1_j)=l_1\delta^{(n)}_{ij}$
if $j\le L_1$. Then we have
\[
E(p)=E(p^1)+E(p^2).
\]
\end{lemma}

\Proof
It is clear by Proposition \ref{prop:E-line}.
\qed

\begin{lemma} \label{lem:3}
Let $\mu$ be any partition. Let $p$ (resp. $\pbar$) be a path
(resp. the ground state path) corresponding to $\mu$. Then we have
\[
E(p)\ge E(\pbar).
\]
\end{lemma}

\Proof
Decompose $\mu$ into $\mu=\mu^1+\cd+\mu^s$ such that each $\mu^J$
is of rectangular shape of distinct depth. From Lemma \ref{lem:1},
we get
\[
E(p)\ge\min_{(p^1,\cd,p^s)}\sum_J E(p^J),
\]
where the minimum is taken over all $(\mu^1,\cd,\mu^s)$-partitions of
$p$. Let $\pbar^J$ be the ground state path corresponding to $\mu^J$.
It is known in the theory of perfect crystals that $E(p^J)\ge E(\pbar^J)$
for any path $p^J$ corresponding to $\mu^J$. Combining with Lemma
\ref{lem:2}, we get
\[
E(p)\ge\sum_J E(\pbar^J)=E(\pbar).
\]
\qed

\begin{lemma} \label{lem:4}
Let $L$ be a positive integer and fix $r$ to be $0\le r<n$. Let
$p=\cd\ot b_2\ot b_1$ (resp. $\pbar=\cd\ot\bbar_2\ot\bbar_1$)
be a path (resp. the ground state path) corresponding to the
partition $(1^L)$. When $L\to\infty$ with $L\equiv r$ mod $n$,
$E(p)-E(\pbar)\to \infty$ unless the set $Y(p)=\{j\mid b_j\neq\bbar_j\}$
is bounded.
\end{lemma}

\Proof
Recall that if a path $p$ corresponding to a single column
$(1^L)$, the energy $E(p)$ is given by $E(p)=\sum_{j=1}^{L-1}
j\theta(b_{j+1}-b_j)$. Here $\theta(z)=1$ ($z\ge0$), $=0$ ($z<0$),
and $b_j$ is identified with a number from $\{1,\cd,n\}$. Fix
$N_0$ such that $b_{N_0n+r+1-j}\neq\bbar_{N_0n+r+1-j}$ for
some $j$ satisfying $1\le j\le n$, and set $L_0=N_0n+r$.
Then there exists $j_0$ such that $1\le j_0<n,
\theta(b_{L_0+1-j_0}-b_{L_0-j_0})=1$. It is easy to see that
in this case
\begin{eqnarray*}
E(p)-E(\pbar)&\ge&\sum_{k=0}^\infty([L_0-kn-j_0]_+-[L_0-kn-n]_+)\\
&\ge&(N_0-1)(n-j_0).
\end{eqnarray*}
Here $[z]_+=z$ ($z\ge0$), $=0$ ($z<0$). Since $N_0$ can be arbitrarily
large, $E(p)-E(\pbar)\to\infty$.
\qed

\begin{lemma} \label{lem:5}
Let $l,N$ be positive integers and $\mu$ be a partition such that
$l\ge\mu_1$. Set $\mu^{(N)}=(l^{Nn})\cup\mu$. Let $p=\cd\ot b_2\ot b_1$
(resp. $\pbar=\cd\ot\bbar_2\ot\bbar_1$) be a path (resp. the ground
state path) corresponding to $\mu^{(N)}$. When $N\to\infty$,
$E(p)-E(\pbar)$ remains finite if and only if the set $Y(p)=\{j\mid
b_j\neq\bbar_j\}$ is bounded.
\end{lemma}

\Proof
The ``if'' part is trivial. To show the ``only if'' part, we decompose
$\mu^{(N)}$ into $\mu^{(N)}=\mu^1+\cd+\mu^s\,
(s=\mu^{(N)}_1)$ such that each $\mu^J$ is of shape $(1^{L_J})$. Let
$(p^1,\cd,p^s)$ be a $(\mu^1,\cd,\mu^s)$-partition of $p$ and $\pbar^J$
be the ground state path corresponding to $\mu^J$. Assume now that
$J(p)$ is not bounded. The previous lemma and inequality
\[
E(p)-E(\pbar)\ge\min_{(p^1,\cd,p^s)}\sum_J(E(p^J)-E(\pbar^J))
\]
shows $E(p)-E(\pbar)\to\infty$, which completes the proof.
\qed

\bigskip
\noindent{\sl Proof of Theorem \ref{th:main_th}.}\quad
We only give proofs when $s=2$, since the
case when $s>2$ is essentially the same.

Let $L_1,L_2$ be finite at first. From Lemma \ref{lem:1} and
\ref{lem:2}, we have an inequality
\[
E(p)-E(\pbar)\ge\min_{(p^1,p^2)}\sum_{J=1}^2(E(p^J)-E(\pbar^J)).
\]
Here the minimum is taken over all
$((l_1^{L_1})\cup\mu^1,(l_2^{L_2})\cup\mu^2)$-partitions of $p$.
Set $p^J=\cd\ot b^J_2\ot b^J_1$ for $J=1,2$.
Now consider the limit when both $L_1-L_2$ and $L_2$ go to infinity,
and let $p$ be an element of $\P_{<\infty}((l_1,r_1,\mu^1),(l_2,r_2,\mu^2))$.
If there is no partition $(p^1,p^2)$ of $p$ such that
$Y(p^J)=\{j\mid b^J_j\neq\bbar^J_j\}$ are bounded both for
$J=1,2$, we have $E(p)-E(\pbar)\to\infty$ by Lemma \ref{lem:5}.
Thus there exists a partition $(p^1,p^2)$ of $p$ such that
$E(p^J)-E(\pbar^J)$ remains finite for $J=1,2$. Since $Y(p^1)$ is
bounded, this $(p^1,p^2)$ is unique in the limit. This shows (1).

(2) is clear from Lemma \ref{lem:2}.
To see (3) note that from Lemma \ref{lem:5}
we have ${\cal P}_{< \infty}((l,r,\mu)) =
{\cal P}_{< \infty}((l,r,\emptyset)) \otimes B_{(\mu_1)}
\otimes \cdots  \otimes B_{(\mu_m)} \, (m = l(\mu))$.
{}From the same lemma and
the theory of perfect crystals \cite{KMN1},
we also have ${\cal P}_{< \infty}((l,r,\emptyset)) \simeq B(l\La^{cl}_r)$,
which proves (3) by the definition (\ref{eq:decomp}).
\qed

\subsection{Affine weight} \label{subsec:aff-wt}
So far we only considered the weight of crystal in the set $\Pcl=
P/\Z\delta$. In this subsection, we take the degree of $\delta$
into account by using the energy of path. This is of crucial
importance in the following sections.

By Theorem \ref{th:main_th}, we reduce our consideration to the set
$\P_{<\infty}((l,r,\mu))$. Let $\pbar$ be the ground state path of
$\P_{<\infty}((l,r,\mu))$. In view of Proposition \ref{prop:interpre},
it is natural to define the {\em affine} weight of $p\in\P_{<\infty}
((l,r,\mu))$ by
\[
\mbox{\sl affine\,wt}\,p=af(\wt p)-(E(p)-E(\pbar))\delta\in P.
\]
For the map $af$, see Section \ref{subsec:preli}. We want
to know the difference of the degree of $\delta$ of each connected
component of $\P_{<\infty}((l,r,\mu))$. Since any highest weight
element $p$ in $\P_{<\infty}((l,r,\mu))$ is of the form
\begin{eqnarray*}
&&p=\pbar^\infty\ot p^\mu,\\
&&\pbar^\infty\mbox{is the ground state path of }
\lim_{N\to\infty}B^{\ot(Nn+r)},\\
&&p^\mu\in{\cal H}(l\La^{cl}_r,\mu),
\end{eqnarray*}
we have
\[
E(p)-E(\pbar)=E(p^\mu)-E(\pbar^\mu),
\]
where $\pbar^\mu\in B_{(\mu_1)}\ot\cd\ot B_{(\mu_m)}$ is defined
{}from $\pbar=\pbar^\infty\ot\pbar^\mu$.
Here we note that the first component of $p^\mu$ is fixed from
the condition for $p^\mu$ to be in ${\cal H}(l\La^{cl}_r,\mu)$.
By abuse of notation we set ${\cal H}(af(\la),\mu)={\cal H}(\la,\mu)$
for $\la\in\Pcl^+$. Calculating $E(\pbar^\mu)$ explicitly, we get

\begin{proposition} \label{prop:P-crystal}
As a $P$-weighted crystal, $\P_{<\infty}((l,r,\mu))$ is isomorphic to
\[
\bigoplus_{p\in{\cal H}(l\La_r,\mu)}
B(l\La_r+af(\wt p)-(E(p)-\ol{E}(l\La_r,\mu))\delta),
\]
where
\begin{equation}
\ol{E}(l\La_r,\mu)=\frac12\sum_{j=1}^{\mu_1}
\left(\frac{t_j^2}{n}-t_j+(\La_{t_j}|\La_{t_j})\right),
\qquad t_j=\mu'_j+r. \label{gsenergy}
\end{equation}
\end{proposition}

\section{One dimensional sums}\label{sect3}

In \cite{KMOTU3}, three kinds of paths,
unrestricted, classically restricted and restricted ones, are
studied, which are associated with the
tensor power of a single perfect crystal.
Here we shall introduce their inhomogeneous versions
associated to $B_{(\mu_1)} \otimes \cdots \otimes B_{(\mu_m)}$
and $B_{(1^{\mu_1})} \otimes \cdots \otimes B_{(1^{\mu_m})}$, and
state the consequence of Theorem \ref{th:main_th}.

\subsection{Unrestricted, classically restricted and restricted paths}
\label{subsec:paths}

Let $\overline{\Lambda}_i = \Lambda_i - \Lambda_0$ be the classical
part of $\Lambda_i$.
We set
$\overline{P} = \oplus_{i=1}^{n-1} {\bf Z} \overline{\Lambda}_i$
and $\overline{P}^+
= \oplus_{i=1}^{n-1} {\bf Z}_{\ge 0} \overline{\Lambda}_i$
as in Section \ref{subsec:preli}.
For any partition $\mu =(\mu_1, \ldots, \mu_m (> 0))$
and $\lambda \in \overline{P}$, we define
\begin{eqnarray}\label{3.1}
{\cal P}_\mu(\lambda) & = & \{
p = b_1 \otimes \cdots \otimes b_m \in
B_{(\mu_1)} \otimes \cdots \otimes B_{(\mu_m)} \mid
af(\wt p)  = \lambda \}, \nonumber \\
{\cal P}'_\mu(\lambda) & = & \{
p = b_1 \otimes \cdots \otimes b_m \in
B_{(1^{\mu_1})} \otimes \cdots \otimes B_{(1^{\mu_m})} \mid
af(\wt p) = \lambda \}.
\end{eqnarray}
Elements of ${\cal P}_\mu(\lambda)$ and
${\cal P}'_\mu(\lambda)$ will be called the
unrestricted paths of weight $\lambda$.
For the latter we assume $\mu_1 \le n-1$.
For $\lambda \in \overline{P}^+ \subset \overline{P}$, we set
\begin{equation}\label{3.2}
\overline{\cal P}_\mu({\lambda})  =  \left\{
b_1 \otimes \cdots \otimes b_m \in
{\cal P}_\mu({\lambda}) \left|
\begin{array}{l}
\varepsilon_i(b_j) \le \langle h_i, \wt b_1 + \cdots + \wt b_{j-1} \rangle  \\
{\rm for }\, 1 \le i \le n-1 \,{\rm and }\,  j = 1, \ldots, m
\end{array} \right. \right\}.
\end{equation}
We define $\overline{\cal P}'_\mu(\lambda)$ similarly by the above
equation by replacing ${\cal P}_\mu({\lambda})$ with
${\cal P}'_\mu({\lambda})$.
Elements of
$\overline{\cal P}_\mu({\lambda})$ and $\overline{\cal P}'_\mu({\lambda})$
will be called the classically restricted paths of weight $\lambda$.
Finally for $\lambda \in P^+_l$  we set
\begin{equation}\label{3.3}
{\cal P}^{(l)}_\mu(\lambda) = \left\{
b_1 \otimes \cdots \otimes b_m \in
{\cal P}_\mu(\overline{\lambda}) \left|
\begin{array}{l}
\varepsilon_i(b_j) \le \langle h_i,
l\Lambda^{cl}_0 + \wt b_1 + \cdots + \wt b_{j-1} \rangle  \\
{\rm for }\, 0 \le i \le n-1 \,{\rm and }\,  j = 1, \ldots, m
\end{array} \right. \right\}.
\end{equation}
We define ${\cal P}^{(l)'}_\mu(\lambda)$ similarly by the above equation
by replacing ${\cal P}_\mu(\overline{\lambda})$ with
${\cal P}'_\mu(\overline{\lambda})$.
Elements of
${\cal P}^{(l)}_\mu(\lambda)$ and ${\cal P}^{(l)'}_\mu(\lambda)$
will be called the (level $l$) restricted paths of weight $\lambda$.
For the former we assume $\mu_1 \le l$.
Notice that
$ {\cal H}(l\Lambda_0,\mu)
= \bigsqcup_{\la\in P^+_l} {\cal P}^{(l)}_\mu(\lambda)$.
See (\ref{eq:defH}).
For $\lambda \in \overline{P}^+$ such that
$\lambda + l \Lambda_0 \in P^+_l$ one has the relations
\begin{eqnarray*}
&&{\cal P}^{(l)}_\mu(\lambda + l \Lambda_0) \subset
\overline{\cal P}_\mu(\lambda) \subset
{\cal P}_\mu(\lambda),\nonumber \\
&&{\cal P}^{(l)'}_\mu(\lambda + l \Lambda_0) \subset
\overline{\cal P}'_\mu(\lambda) \subset
{\cal P}'_\mu(\lambda).
\end{eqnarray*}

Having defined the paths we now introduce the
associated one dimensional sums (1dsums) by
\begin{eqnarray}\label{3.4}
g_\mu(\lambda) &= \sum_{p \in {\cal P}_\mu(\lambda)} q^{E(p)}
\quad & \lambda \in \overline{P},\nonumber\\
g'_\mu(\lambda) &= \sum_{p \in {\cal P}'_\mu(\lambda)} q^{-E(p)}
\quad & \lambda \in \overline{P},\nonumber\\
{X}_\mu({\lambda})&=
\sum_{p \in \overline{\cal P}_\mu({\lambda})} q^{E(p)}
\quad & {\lambda} \in \overline{P}^+,\nonumber\\
{X'}_\mu({\lambda})&=
\sum_{p \in \overline{\cal P}'_\mu({\lambda})} q^{-E(p)}
\quad & {\lambda} \in \overline{P}^+,\nonumber\\
{X}^{(l)}_\mu(\lambda)&=
\sum_{p \in {\cal P}^{(l)}_\mu(\lambda)} q^{E(p)}
\quad & \lambda \in P^+_l,\nonumber \\
{X}^{(l)'}_\mu(\lambda)&=
\sum_{p \in {\cal P}^{(l)'}_\mu(\lambda)} q^{-E(p)}
\quad & \lambda \in P^+_l,
\end{eqnarray}
where the energy $E(p)$ of a path $p$ is given in
Section \ref{subsec:energy}.
We note that when $\forall \mu_i = 1$, the energy functions in the two
cases are different by an additive constant.
The functions $g_\mu(\lambda),
{X}_\mu({\lambda})$,
${X}^{(l)}_\mu(\lambda)$ (and their analogues for
$B_{(1^{\mu_1})} \otimes \cdots \otimes B_{(1^{\mu_m})}$) are called
the unrestricted, the classically restricted and the
(level $l$) restricted 1dsum,
respectively.
They are polynomials in $q$ with non-negative integer
coefficients.
For rectangular $\mu$, the unrestricted and
the restricted 1dsums have the origin in the studies of
solvable vertex and RSOS models in statistical mechanics
\cite{ABF}, \cite{DJKMO1}, \cite{DJKMO2}, \cite{JMO}.
In \cite{KMOTU3}, paths associated with
$B \otimes \cdots \otimes B$ were studied from the
Demazure module viewpoint.
They are homogeneous but $B$ can be
{\em any} perfect crystal.
Let $g_j(b,\la)_B, \overline{X}_j(b, \xi, \eta)_B$ and
$X_j(b, \xi, \eta)_B$ be the unrestricted, classically restricted
and restricted 1dsums introduced there, respectively,
where we have exhibited the $B$-dependence explicitly.
They are related with the 1dsums in (\ref{3.4}) as
\begin{eqnarray}
&g_{(l^j)}(\la)  = q^{-jl}g_j(b^+,\la)_{B_{(l)}},
&g'_{(l^j)}(\la)  = g_j(b^-,\la)_{B_{(1^l)}},
\nonumber \\
&X_{(l^j)}(\la)  =
q^{-jl}\overline{X}_j(b^+,0,\la)_{B_{(l)}},
&X'_{(l^j)}(\la)  =
\overline{X}_j(b^-,0,\la)_{B_{(1^l)}},
\nonumber \\
&X^{(l+l')}_{(l^j)}(\la)  =
q^{-jl}X_j(b^+,l\La_0,\la)_{B_{(l)}}, \quad
&X^{(l+l')'}_{(l^j)}(\la)  =
X_j(b^-,l\La_0,\la)_{B_{(1^l)}},
\nonumber
\end{eqnarray}
where
$b^+ = (\overbrace{0, \ldots, 0}^{n-1},l) \in B_{(l)}$
and
$b^- =
(\overbrace{0, \ldots, 0}^{n-l},
\overbrace{1, \ldots, 1}^{l}) \in B_{(1^l)}$.

% ************************* Example of paths (symmetric case) *************
\noindent
\begin{example}\label{ex1path}
Let $\lambda = (321), \; \mu = (2211)$ and $n=3$.
All the elements $p$ of
${\cal P}_\mu(\lambda),\; \overline{\cal P}_\mu(\lambda)$ and
${\cal P}^{(2)}_\mu(\lambda)$ are listed with their energy $E(p)$.
The elements of
$\overline{\cal P}_\mu(\lambda)$ and ${\cal P}^{(2)}_\mu(\lambda)$
are labelled c and r, respectively.
The element $(x_1,x_2,x_3)$ of each crystal $B_{(k)}$ is denoted by
$1^{x_1}2^{x_2}3^{x_3}$.
  \begin{center}
    \begin{tabular}{|c|c|c||c|c|c|}   \hline
      $p$ & $E$ & \hspace{2,5ex} & $p$ & $E$ & \hspace{2,5ex} \\ \hline \hline
      $ 11 \otimes 12 \otimes 2\otimes 3 $ & $3$ & c & %
      $ 11 \otimes 12 \otimes 3\otimes 2 $ & $2$ & c \\
      $ 11 \otimes 13 \otimes 2\otimes 2 $ & $4$ &   & %
      $ 11 \otimes 22 \otimes 1\otimes 3 $ & $2$ & c \\
      $ 11 \otimes 22 \otimes 3\otimes 1 $ & $1$ & c,r & %
      $ 11 \otimes 23 \otimes 1\otimes 2 $ & $2$ &   \\
      $ 11 \otimes 23 \otimes 2\otimes 1 $ & $3$ &   & %
      $ 12 \otimes 11 \otimes 2\otimes 3 $ & $4$ &   \\
      $ 12 \otimes 11 \otimes 3\otimes 2 $ & $5$ &   & %
      $ 12 \otimes 12 \otimes 1\otimes 3 $ & $3$ &   \\
      $ 12 \otimes 12 \otimes 3\otimes 1 $ & $2$ &   & %
      $ 12 \otimes 13 \otimes 1\otimes 2 $ & $3$ &   \\
      $ 12 \otimes 13 \otimes 2\otimes 1 $ & $4$ &   & %
      $ 12 \otimes 23 \otimes 1\otimes 1 $ & $3$ &   \\
      $ 13 \otimes 11 \otimes 2\otimes 2 $ & $5$ &   & %
      $ 13 \otimes 12 \otimes 1\otimes 2 $ & $3$ &   \\
      $ 13 \otimes 12 \otimes 2\otimes 1 $ & $4$ &   & %
      $ 13 \otimes 22 \otimes 1\otimes 1 $ & $4$ &   \\
      $ 22 \otimes 11 \otimes 1\otimes 3 $ & $6$ &   & %
      $ 22 \otimes 11 \otimes 3\otimes 1 $ & $5$ &   \\
      $ 22 \otimes 13 \otimes 1\otimes 1 $ & $4$ &   & %
      $ 23 \otimes 11 \otimes 1\otimes 2 $ & $6$ &   \\
      $ 23 \otimes 11 \otimes 2\otimes 1 $ & $5$ &   & %
      $ 23 \otimes 12 \otimes 1\otimes 1 $ & $7$ &   \\ \hline
    \end{tabular}
  \end{center}
\end{example}
%
% ************************ Example of paths (anti-symmetric case) *************
\vspace{5mm}
\begin{example}\label{ex2path}
Let $\lambda = (322), \; \mu = (2221)$ and $n=3$.
All the elements $p$ of
${\cal P}'_\mu(\lambda),\; \overline{\cal P}'_\mu(\lambda)$ and
${\cal P}^{(1)'}_\mu(\lambda)$ are listed with their energy $E(p)$.
The elements of $\overline{\cal P}'_\mu(\lambda)$ and
${\cal P}^{(1)'}_\mu(\lambda)$
are labelled c and r, respectively.
The element $(x_1,x_2,x_3)$ of each crystal $B_{(1^k)}$ is denoted by
$1^{x_1}2^{x_2}3^{x_3}$.
  \begin{center}
    \begin{tabular}{|c|c|c||c|c|c|}   \hline
      $p$ & $E$ & \hspace{2,5ex} & $p$ & $E$ & \hspace{2,5ex} \\ \hline \hline
      $ 12 \otimes 12 \otimes 13 \otimes 3 $ & $-2$ & c & %
      $ 12 \otimes 13 \otimes 12 \otimes 3 $ & $-3$ & c \\
      $ 12 \otimes 13 \otimes 13 \otimes 2 $ & $-3$ &   & %
      $ 12 \otimes 13 \otimes 23 \otimes 1 $ & $-4$ & c,r \\
      $ 12 \otimes 23 \otimes 13 \otimes 1 $ & $-2$ &   & %
      $ 13 \otimes 12 \otimes 12 \otimes 3 $ & $-1$ &   \\
      $ 13 \otimes 12 \otimes 13 \otimes 2 $ & $-2$ &   & %
      $ 13 \otimes 12 \otimes 23 \otimes 1 $ & $-1$ &   \\
      $ 13 \otimes 13 \otimes 12 \otimes 2 $ &  $0$ &   & %
      $ 13 \otimes 23 \otimes 12 \otimes 1 $ & $-2$ &   \\
      $ 23 \otimes 12 \otimes 13 \otimes 1 $ & $-1$ &   & %
      $ 23 \otimes 13 \otimes 12 \otimes 1 $ &  $0$ &   \\ \hline
    \end{tabular}
  \end{center}
\end{example}

\subsection{1dsums and Kostka-Foulkes polynomials} \label{subsec:1dsum}

Given a partition $\mu$ and
$\lambda \in P$, it is easy to see that the sets (\ref{3.1})--(\ref{3.3})
are empty hence the 1dsums (\ref{3.4}) are zero, unless
$\vert \mu \vert -
\langle \sum_{j=1}^{n-1} j h_j,  \lambda \rangle
\in n{\bf Z}_{\ge 0}$.
In view of this we shall identify
$\lambda \in P$ with the composition
$(\lambda_1, \ldots, \lambda_n) \in ({\bf Z}_{\ge 0})^n$ such that
$\overline{\la} = \sum_{i=1}^{n-1}(\la_i-\la_{i+1})\overline{\Lambda}_i$
and $\sum_{i=1}^n \lambda_i = \vert \mu \vert$.
This is independent of the level of $\lambda$.
If $\lambda \in \overline{P}^+$, the composition
$(\lambda_1, \ldots, \lambda_n)$ is in fact a partition;
$\lambda_1 \ge \cdots \ge \lambda_n$.
If $\lambda \in P^+_l$, it is a (level $l$) {\it restricted}
partition, by which we mean
$l + \lambda_n \ge \lambda_1 \ge \cdots \ge \lambda_n$.
We regard a partition
$\la = (\la_1, \ldots, \la_n)$ as a Young diagram
in which the length of the $i$th row is $\la_i$.
The depth of $\la$ will be denoted by $l(\la)$.
The conjugate partition $\la'$ is obtained from $\la$ by
transposing the diagram.
Given partitions $\xi$ and $\eta$ let
$K_{\xi \eta}(q)$ be the associated Kostka-Foulkes polynomial and
$K_{\xi \eta} = K_{\xi \eta}(1)$ the Kostka number \cite{Ma}.
We extend the definition of the latter to
compositions $\eta \in ({\bf Z}_{\ge 0})^n$ by assuming the
invariance under any permutation of
the components of $\eta$.
Under the above identification of
$\lambda \in P$ with the compositions the 1dsums are expressed as
\begin{proposition} \label{pr:1dsumbyk}
\begin{eqnarray}
g_\mu(\lambda) & = &\sum_\eta K_{\eta \lambda}
K_{\eta \mu}(q) \quad \quad \lambda \in \overline{P},
\label{gbyk}\\
g'_\mu(\lambda) & = &\sum_\eta K_{\eta' \lambda}
K_{\eta \mu}(q) \quad \quad \lambda \in \overline{P},
\label{gbykp}\\
{X}_\mu(\lambda) & = & K_{\lambda \mu}(q), \quad
{X}'_\mu(\lambda)  =  K_{\lambda' \mu}(q)
\quad\lambda \in \overline{P}^+,\label{xbarbyk}
\end{eqnarray}
where the sum in (\ref{gbyk}) (resp. (\ref{gbykp})) runs
over all partitions $\eta$ of
$\vert \mu \vert $ with $\l(\eta) \le n$ (resp. $\eta_1 \le n$).
\end{proposition}

\Proof
(\ref{xbarbyk}) is due to \cite{NY}.
To see (\ref{gbyk})
consider the character of the energy and the
$\goth{sl}_{\, n}$ weight over
$B_{(\mu_1)} \otimes \cdots \otimes B_{(\mu_m)}$.
It can be counted either weight vector-wise or
$\goth{sl}_{\, n}$ component-wise,
leading to
$\sum_{\lambda}g_\mu(\lambda) m_\lambda =
\sum_\eta K_{\eta \mu}(q) s_\eta$.
Here $m_\la$ and $s_\eta$ denote the monomial symmetric function and
the Schur function \cite{Ma}, respectively
and the both sums run over the partitions of $\vert \mu \vert$
with depth $\le n$.
Substituting the expansion $s_\eta = \sum_\la K_{\eta \la}m_\la$
into this and comparing the coefficients of $m_\la$ one has (\ref{gbyk}).
(\ref{gbykp}) can be verified similarly.
\qed

% ********** Example of KK proposition (symmetric case) ********************
\noindent
\begin{example}\label{ex1kk}
Let $\lambda = (321), \; \mu = (2211)$ and $n=3$ as in
Example \ref{ex1path}.
The RHS of (\ref{gbyk}) is calculated from the data
\begin{center}
\begin{tabular}{|c|c|c|}
\hline
$\eta$ & $K_{\eta (321)}$ & $K_{\eta (2211)}(q)$ \\
\hline
$(6)   $ & $1$ & $q^7         $ \\
$(51)  $ & $2$ & $q^4+q^5+q^6 $ \\
$(42)  $ & $2$ & $2q^3+q^4+q^5$ \\
$(41^2)$ & $1$ & $q^2+q^3+q^4 $ \\
$(3^2) $ & $1$ & $q^2+q^4     $ \\
$(321) $ & $1$ & $q+2q^2+q^3  $ \\
\hline
\end{tabular}
\end{center}
as
\begin{eqnarray}
\mbox{RHS} &=& \sum_{\eta \; (l(\eta) \le 3)}
K_{\eta (321)} K_{\eta (2211)}(q)\nonumber\\
&=& q + 4q^2 + 6q^3 + 6q^4 + 4q^5 + 2q^6 + q^7,\label{qpol1}
\end{eqnarray}
in agreement with Example \ref{ex1path}.
\end{example}

% ********** Example of KK proposition (anti-symmetric case) ***************
\noindent
\begin{example}\label{ex2kk}
Let $\lambda = (322), \; \mu = (2221)$ and $n=3$ as in Example \ref{ex2path}.
The RHS of (\ref{gbykp}) is calculated from the data
\begin{center}
\begin{tabular}{|c|c|c|}
\hline
$\eta$ & $K_{\eta' (322)}$ & $K_{\eta (2221)}(q)$ \\
\hline
$(3^21) $ & $1$ & $q^2+q^3+q^4$ \\
$(32^2) $ & $1$ & $q+q^2+q^3  $ \\
$(321^2)$ & $2$ & $q+q^2      $ \\
$(2^31) $ & $2$ & $1          $ \\
\hline
\end{tabular}
\end{center}
as
\begin{eqnarray}
\mbox{RHS} &=& \sum_{\eta \; (\eta_1 \le 3)}
K_{\eta' (322)} K_{\eta (2221)}(q)\nonumber\\
&=& 2 + 3q + 4q^2 + 2q^3 + q^4, \label{qpol2}
\end{eqnarray}
in agreement with Example \ref{ex2path}.
\end{example}

By analogy with (\ref{xbarbyk}) the
1dsums $X^{(l)}_\mu(\la)$ and  $X^{(l)'}_\mu(\la)$ are also
called the (level $l$) restricted Kostka-Foulkes polynomials.
%\begin{equation}\label{rkostka}
%K^{(l)}_{\la \mu}(q) = X^{(l)}_\mu(\la), \quad
%K^{(l)}_{\la' \mu}(q) = X^{(l)'}_\mu(\la),
%\end{equation}
%where $\la$ is any level $l$ restricted partition of $\vert \mu \vert$.

\subsection{Limit of 1dsums} \label{subsec:1dsumlimit}
Given an $U_q(\widehat{\goth{sl}}_{\, n})$-module $M$ and $\la \in P$,
we define
\begin{eqnarray*}
M_\la & = & \{ v \in M \mid \wt v = \la \},\\
\mbox{$[ M : \la ]_{cl}$} & = & \mbox{dim}
\{ v \in M_\la \mid e_i v = 0 \hbox{ for } 1 \le i \le n-1 \},\\
\mbox{$[ M : \la ]$} & = & \mbox{dim}
\{ v \in M_\la \mid e_i v = 0 \hbox{ for } 0 \le i \le n-1 \}.
\end{eqnarray*}
Suppose that $M$ is a level $l$ module (not necessarily irreducible).
We prepare the following notations
for the relevant branching functions:
\begin{eqnarray}
c^M_\la(q) & = &
\sum_i \left( \dim M_{\la + l\La_0 - i \delta} \right) q^i
\ \ \qquad\qquad \la \in \overline{P},
\label{defc}\\
b^M_\la(q) & = & \sum_i [M :
\la + l\Lambda_0 - i \delta]_{cl}\, q^i
\, \  \quad\quad\quad \la \in \overline{P}^+,
\label{defb}\\
a^M_\la(q) & = &
\sum_i [M : \la - i \delta] \, q^i
\qquad\quad\qquad\qquad \la \in P^+_l.
\label{defa}
\end{eqnarray}
Up to an overall power of $q$, (\ref{defc}) is the
string function \cite{KP},
(\ref{defa}) is a branching coefficient of the module
$V(\la)$ within $M$, and
(\ref{defb}) is the
branching coefficient of
the irreducible $U_q(\goth{sl}_{\, n})$
module with highest weight $\la \in \overline{P}^+$ within $M$,
where $U_q(\goth{sl}_{\, n})$ stands for the subalgebra of
$U_q(\widehat{\goth{sl}}_{\, n})$ generated by $e_i, f_i, t_i
( 1 \le i \le n-1)$.

The $U_q(\widehat{\goth{sl}}_{\, n})$-module corresponding to
Theorem \ref{th:main_th} (1) is given by
\begin{equation}\label{vcal}
{\cal V} = \bigotimes_{J=1}^s \Bigl(
\bigoplus_{p \in {\cal H}(l_J\Lambda_{r_J}, \mu^{J})}
V\bigl(l_J\Lambda_{r_J} + af(\wt p) -
(E(p) - \overline{E}(l_J\Lambda_{r_J}, \mu^{J}))\delta \bigr) \Bigr).
\end{equation}
As a consequence of Theorem \ref{th:main_th} and Proposition \ref{prop:P-crystal},
the large $\mu$ limit (\ref{eq:mu}) of the 1dsums
$g_\mu(\la), {X}_\mu(\la)$ and
$X^{(l)}_\mu(\la)$ gives rise to the branching functions
related to ${\cal V}$.
We summarize them in
\begin{proposition}\label{1dsumlimit}
Let $l_0 = l - \sum_{J=1}^s l_J$ and assume that
$l_0 \ge 0$. Then we have
\begin{eqnarray}
\lim_\mu q^{-\overline{E}(l\La_0,\mu)} g_\mu(\la) & = &
c^{\cal V}_\la(q)
\ \quad\qquad\qquad\la \in \overline{P},
\label{glimit}\\
\lim_{\mu} q^{-\overline{E}(l\La_0,\mu)}
{X}_{\mu}(\la) & = & b^{\cal V}_\la(q)
\ \quad\qquad\qquad \la \in \overline{P}^+,
\label{xbarlimit}\\
\lim_{\mu} q^{-\overline{E}(l\La_0,\mu)}
X^{(l)}_\mu(\lambda) & = &
a^{{\cal V} \otimes V(l_0\Lambda_0)}_\la(q)
\quad\quad \la \in P^+_l.
\label{xlimit}
\end{eqnarray}
\end{proposition}
The $U_q(\widehat{\goth{sl}}_{\, n})$-module ${\cal V}$
(\ref{vcal}) is reducible in general and has the form
\begin{equation}
{\cal V} = \bigoplus \left(V(\xi_1 - d_1 \delta) \otimes \cdots
\otimes V(\xi_s - d_s \delta) \right), \label{vcal2}
\end{equation}
where the direct sum runs over the $s$-tuple of
the restricted paths
$(p^1, \ldots, p^s) \in
{\cal H}(l_1\La_{r_1}, \mu^1) \times \cdots \times
{\cal H}(l_s\La_{r_s}, \mu^s)$.
$\xi_J \in P^+_{l_J}$ and
$d_J \in {\bf Z}$ are specified as
$\xi_J = l_J\La_{r_J} + af(\wt p^J)$ and
$d_J = E(p^J) - \overline{E}(l_J\La_{r_J},\mu^J)$
for each summand.
Correspondingly the branching functions appearing
in (\ref{glimit})--(\ref{xlimit}) are in fact
equal to the linear combinations:
\begin{eqnarray}
c^{\cal V}_\la(q) & = &
\sum q^{d_1 + \cdots + d_s}
c^{V(\xi_1)\otimes \cdots \otimes V(\xi_s)}_\la(q),
\label{glimit2} \\
b^{\cal V}_\la(q) & = &
\sum q^{d_1 + \cdots + d_s}
b^{V(\xi_1)\otimes \cdots \otimes V(\xi_s)}_\la(q),
\label{xbarlimit2} \\
a^{{\cal V}\otimes V(l_0\La_0)}_\la(q) & = &
\sum q^{d_1 + \cdots + d_s}
a^{V(\xi_1)\otimes \cdots \otimes V(\xi_s)\otimes V(l_0\La_0)}_\la(q),
\label{xlimit2}
\end{eqnarray}
where the three sums are taken in the same way as explained above.
In this way our construction by means of the inhomogeneous paths
naturally gives rise to the finite linear combinations of
the ``usual'' branching functions $a, b$ and $c$.
By choosing $l_J, r_J, \mu^J$ ($1 \le J \le s$) variously,
one can let ${\cal V}$ cover a large
family of (generally reducible) $U_q(\widehat{\goth{sl}}_{\, n})$-modules.
In general the corresponding sum
(\ref{glimit2})--(\ref{xlimit2}) is complicated depending
especially on $\mu^1, \ldots, \mu^s$.
We shall therefore calculate the limits
in Proposition \ref{1dsumlimit} explicitly only for $s = 1$ or
$\mu^1 = \cdots = \mu^s = \emptyset$
in Sections \ref{sec:fflimit} and \ref{discussion}.
These sample calculations already lead to
several character formulae, which appear to be new.
To exploit the full variety of such formulae
for general ${\cal V}$ will be an interesting problem.
Now we proceed to
another essential input for our calculation,
{\em fermionic forms} of the 1dsums.

\section{Fermionic formulae}\label{sect4}

Here we shall present explicit formulae for the
1dsums.
They all have a fermionic form, by which we roughly mean the
polynomials or series that are free of signs or possibly have some
relevance to the Bethe ansatz.
Such formulae for $g_\mu(\la)$ and $g'_\mu(\la)$ are
new, while for $X_\mu(\la)$ and $X^{(l)}_\mu(\la)$
they have been known or conjectured in various cases in
\cite{BMS}, \cite{Ki}, \cite{Ki2}, \cite{KR}.
We shall also present
their limiting forms as $\vert \mu \vert \rightarrow \infty$
as in (\ref{eq:mu}).
By virtue of Proposition \ref{1dsumlimit}
this establishes several character formulae concerning
the tensor products of $U_q(\widehat{\goth{sl}}_{\, n})$-modules.
It turns out that the limit of $g_\mu(\la)$ gives rise to a
generalization of the formulae in \cite{FS}, \cite{Ge2}, \cite{KNS}
and those of ${X}_\mu(\la), X^{(l)}_\mu(\la)$ yield
generalizations of a spinon character formulae in
\cite{ANOT}, \cite{BPS}, \cite{BLS} and \cite{NY2}.

In the working below we shall employ the notations:
\begin{eqnarray*}
&\left[ \begin{array}{c} m \\ k \end{array} \right] &=
\left \{\begin{array}{cl}
{\displaystyle (q)_m \over (q)_k (q)_{m-k}} & \quad
\mbox{if \, $0 \le k \le m$}\\
0& \quad \mbox{otherwise}
\end{array}\right., \\
&(q)_m &= \prod_{i=1}^m(1-q^i) \quad {\rm for }\ \  m \in {\bf Z}_{\ge 0}.
\end{eqnarray*}
We shall also use
\[
 n(\nu) = \sum_{i\ge 1}
\left( \begin{array}{c} \nu'_i \\ 2 \end{array} \right), \qquad
\left( \begin{array}{c} x \\ 2 \end{array} \right) = \frac{x(x-1)}{2}
\]
according to \cite{Ma}.
The symbol $n(\nu)$ should not be confused with
$n$ from $\goth{sl}_{\, n}$.
Elements of the Cartan matrix of $\goth{sl}_{\, k}$ and its inverse will be
written as ($1 \le i, j \le k-1$)
\begin{equation}\label{cartan}
C^{(k)}_{i j} = 2\delta_{i j} - \delta_{i j+1} - \delta_{i j-1}, \quad
C^{(k)^{-1}}_{i j} = \mbox{ min}(i,j) - \frac{i j}{k}.
\end{equation}
We assume that
$C^{(k)^{-1}}_{i j} = 0$ if $i$ or $j$ is equal to
$0$ or $k$.

\subsection{Fermionic formulae of 1dsums}\label{ssec4.1}

A fermionic formula of the unrestricted 1dsum $g_\mu(\la)$
(\ref{gbyk}) is given by
\begin{proposition}\label{pro:ffkk}
For any partition $\mu$ and composition
$\la \in ({\bf Z}_{\ge 0})^n$ such that
$\vert \la \vert = \vert \mu \vert$ we have
\begin{eqnarray}\label{eq:ffkk}
\sum_{\eta \, (l(\eta) \le n)} K_{\eta \la}K_{\eta \mu}(q) & = &
\sum_{\{\nu \}} q^{\phi(\{\nu \})}
\prod_{{\scriptstyle 1 \le a \le n-1} \atop
   {\scriptstyle 1 \le i \le \mu_1}}
\left[ \begin{array}{c} \nu^{(a+1)}_i -  \nu^{(a)}_{i+1}
 \\   \nu^{(a)}_i -  \nu^{(a)}_{i+1}
\end{array} \right],\\
\phi(\{\nu \}) & = & \sum_{a=0}^{n-1} \sum_{i=1}^{\mu_1}
\left( \begin{array}{c} \nu^{(a+1)}_i -  \nu^{(a)}_i \\
2 \end{array} \right), \label{eq:ffkkphi}
\end{eqnarray}
where the sum $\sum_{\{\nu \}}$ runs over the sequences of
Young diagrams $\nu^{(1)}, \ldots, \nu^{(n-1)}$ such that
\begin{eqnarray}\label{eq:ffkknu}
&&\emptyset =: \nu^{(0)} \subset \nu^{(1)} \subset \cdots \subset
\nu^{(n-1)} \subset \nu^{(n)}:= \mu', \nonumber \\
&&\vert \nu^{(a)} \vert = \la_1 + \cdots + \la_a \quad
\mbox{ for } 1 \le a \le n-1.\nonumber\\
\end{eqnarray}
\end{proposition}

When $n=2$ and $\mu$ is rectangular
the result
$g_\mu(\la) =$ RHS of (\ref{eq:ffkk}) was known in \cite{DJKMO1},
and the above formula is reduced to the one in \cite{Ki}.
For $n=2$ and general $\mu$
it agrees with  \cite{SW}.
The case $\mu=(1^m)$ was known for general $n$
in \cite{JMO} and \cite{DJKMO2}.
See also (\ref{proofkk})--(\ref{proofkkwa})
for an equivalent expression
in terms of ``TBA-like'' variables and
Remark \ref{rem:tba} concerning it.
%
%A proof of Proposition \ref{pro:ffkk} will be given
%in our forthcoming paper
%\cite{HKKOTY}.
The proof of Proposition \ref{pro:ffkk} will be given
after the proof of Proposition \ref{pro:ffkkp}.
%%%%%%%%%%%%%%%%%%%%%%%%%%yyy
%
Here we shall illustrate the formula by an example.

\noindent
\begin{example}\label{ex1ff}
Let $\lambda = (321), \; \mu = (2211)$ and $n=3$ as in
Example \ref{ex1path} and \ref{ex1kk}.
The relevant $\{ \nu^{(1)}, \nu^{(2)}, \nu^{(3)} = (42) \}$'s
and their contributions are
\\ \noindent
\setlength{\unitlength}{3mm}
% FF sym 1
\begin{picture}(20,3.5)(-10,-0.5)
\put(0,2){\line(1,0){3}}
\put(0,1){\line(1,0){3}}
\put(0,1){\line(0,1){1}}
\put(1,1){\line(0,1){1}}
\put(2,1){\line(0,1){1}}
\put(3,1){\line(0,1){1}}
\put(3,0){\makebox(2,2){$ \subset$}}
\put(5,2){\line(1,0){3}}
\put(5,1){\line(1,0){3}}
\put(5,0){\line(1,0){2}}
\put(5,0){\line(0,1){2}}
\put(6,0){\line(0,1){2}}
\put(7,0){\line(0,1){2}}
\put(8,1){\line(0,1){1}}
\put(9,0){\makebox(2,2){$ \subset$}}
\put(11,2){\line(1,0){4}}
\put(11,1){\line(1,0){4}}
\put(11,0){\line(1,0){2}}
\put(11,0){\line(0,1){2}}
\put(12,0){\line(0,1){2}}
\put(13,0){\line(0,1){2}}
\put(14,1){\line(0,1){1}}
\put(15,1){\line(0,1){1}}
\put(16,0){\makebox(4,2){$\quad q^4 {\displaystyle {2 \brack 1}}$,}}
\end{picture}
% FF sym 2
\par\noindent
\begin{picture}(20,3.5)(-10,-0.5)
\put(0,2){\line(1,0){3}}
\put(0,1){\line(1,0){3}}
\put(0,1){\line(0,1){1}}
\put(1,1){\line(0,1){1}}
\put(2,1){\line(0,1){1}}
\put(3,1){\line(0,1){1}}
\put(3,0){\makebox(2,2){$ \subset$}}
\put(5,2){\line(1,0){4}}
\put(5,1){\line(1,0){4}}
\put(5,0){\line(1,0){1}}
\put(5,0){\line(0,1){2}}
\put(6,0){\line(0,1){2}}
\put(7,1){\line(0,1){1}}
\put(8,1){\line(0,1){1}}
\put(9,1){\line(0,1){1}}
\put(9,0){\makebox(2,2){$ \subset$}}
\put(11,2){\line(1,0){4}}
\put(11,1){\line(1,0){4}}
\put(11,0){\line(1,0){2}}
\put(11,0){\line(0,1){2}}
\put(12,0){\line(0,1){2}}
\put(13,0){\line(0,1){2}}
\put(14,1){\line(0,1){1}}
\put(15,1){\line(0,1){1}}
\put(16,0){\makebox(4,2){$\quad q^3 {\displaystyle {4 \brack 3}{2 \brack 1}}$,}}
\end{picture}
% FF sym 3
\par\noindent
\begin{picture}(20,3.5)(-10,-0.5)
\put(0,2){\line(1,0){2}}
\put(0,1){\line(1,0){2}}
\put(0,0){\line(1,0){1}}
\put(0,0){\line(0,1){2}}
\put(1,0){\line(0,1){2}}
\put(2,1){\line(0,1){1}}
\put(3,0){\makebox(2,2){$ \subset$}}
\put(5,2){\line(1,0){4}}
\put(5,1){\line(1,0){4}}
\put(5,0){\line(1,0){1}}
\put(5,0){\line(0,1){2}}
\put(6,0){\line(0,1){2}}
\put(7,1){\line(0,1){1}}
\put(8,1){\line(0,1){1}}
\put(9,1){\line(0,1){1}}
\put(9,0){\makebox(2,2){$ \subset$}}
\put(11,2){\line(1,0){4}}
\put(11,1){\line(1,0){4}}
\put(11,0){\line(1,0){2}}
\put(11,0){\line(0,1){2}}
\put(12,0){\line(0,1){2}}
\put(13,0){\line(0,1){2}}
\put(14,1){\line(0,1){1}}
\put(15,1){\line(0,1){1}}
\put(16,0){\makebox(4,2)
{$\quad q^2 {\displaystyle {3 \brack 1}{2 \brack 1}}$,}}
\end{picture}
% FF sym 4
\par\noindent
\begin{picture}(20,3.5)(-10,-0.5)
\put(0,2){\line(1,0){2}}
\put(0,1){\line(1,0){2}}
\put(0,0){\line(1,0){1}}
\put(0,0){\line(0,1){2}}
\put(1,0){\line(0,1){2}}
\put(2,1){\line(0,1){1}}
\put(3,0){\makebox(2,2){$ \subset$}}
\put(5,2){\line(1,0){3}}
\put(5,1){\line(1,0){3}}
\put(5,0){\line(1,0){2}}
\put(5,0){\line(0,1){2}}
\put(6,0){\line(0,1){2}}
\put(7,0){\line(0,1){2}}
\put(8,1){\line(0,1){1}}
\put(9,0){\makebox(2,2){$ \subset$}}
\put(11,2){\line(1,0){4}}
\put(11,1){\line(1,0){4}}
\put(11,0){\line(1,0){2}}
\put(11,0){\line(0,1){2}}
\put(12,0){\line(0,1){2}}
\put(13,0){\line(0,1){2}}
\put(14,1){\line(0,1){1}}
\put(15,1){\line(0,1){1}}
\put(16,0){\makebox(4,2){$\quad q {\displaystyle {2 \brack 1}^3}$.}}
\end{picture}
\\
By summing them up, we obtain
(\ref{qpol1}).
\end{example}

A fermionic formula of the unrestricted 1dsum $g'_\mu(\la)$
(\ref{gbykp}) is given by
\begin{proposition}\label{pro:ffkkp}
For any partition $\mu$ and composition
$\la \in ({\bf Z}_{\ge 0})^n$ such that
$\vert \la \vert = \vert \mu \vert$ and
$\mu_1 \le n-1$ we have
\begin{equation}
\sum_{\eta \, (\eta_1 \le n)} K_{\eta' \la}K_{\eta \mu}(q) =
\sum_{\{\nu \}}
\prod_{{\scriptstyle 1 \le a \le n - 1} \atop
   {\scriptstyle 1 \le i \le \mu_1}}
\left[ \begin{array}{c} \nu^{(a+1)}_i -  \nu^{(a+1)}_{i+1}
 \\   \nu^{(a)}_i -  \nu^{(a+1)}_{i+1}
\end{array} \right],\label{eq:ffkkp}
\end{equation}
where the sum $\sum_{\{\nu \}}$ runs over the sequences of
Young diagrams $\nu^{(1)}, \ldots, \nu^{(n-1)}$ such that
\begin{eqnarray}\label{eq:ffkknua}
&&\emptyset =: \nu^{(0)} \subset \nu^{(1)} \subset \cdots \subset
\nu^{(n-1)} \subset \nu^{(n)}:= \mu', \nonumber \\
&&\nu^{(a)}/\nu^{(a-1)}: \quad \mbox{ horizontal strip of length } \la_a.
\nonumber\\
\end{eqnarray}
\end{proposition}
%
%Again we postpone a proof of Proposition \ref{pro:ffkkp}
%to \cite{HKKOTY} and give here an example instead.

\noindent
% ********** Example of FF (anti-symmetric case) **************************
\noindent
\begin{example}\label{ex2ff}
Let $\lambda = (322),\; \mu = (2221)$ and $n=3$ as in
Example \ref{ex2path} and \ref{ex2kk}.
The relevant $\{ \nu^{(1)}, \nu^{(2)}, \nu^{(3)} = (43) \}$'s
and their contributions are
\\ \noindent
\setlength{\unitlength}{3mm}
% FF anti-sym 1
\begin{picture}(20,3.5)(-10,-0.5)
\put(0,2){\line(1,0){3}}
\put(0,1){\line(1,0){3}}
\put(0,1){\line(0,1){1}}
\put(1,1){\line(0,1){1}}
\put(2,1){\line(0,1){1}}
\put(3,1){\line(0,1){1}}
\put(3,0){\makebox(2,2){$ \subset$}}
\put(5,2){\line(1,0){3}}
\put(5,1){\line(1,0){3}}
\put(5,0){\line(1,0){2}}
\put(5,0){\line(0,1){2}}
\put(6,0){\line(0,1){2}}
\put(7,0){\line(0,1){2}}
\put(8,1){\line(0,1){1}}
\put(9,0){\makebox(2,2){$ \subset$}}
\put(11,2){\line(1,0){4}}
\put(11,1){\line(1,0){4}}
\put(11,0){\line(1,0){3}}
\put(11,0){\line(0,1){2}}
\put(12,0){\line(0,1){2}}
\put(13,0){\line(0,1){2}}
\put(14,0){\line(0,1){2}}
\put(15,1){\line(0,1){1}}
\put(16,0){\makebox(4,2){$\quad {\displaystyle {3 \brack 2}}$,}}
\end{picture}
% FF anti-sym 2
\par\noindent
\begin{picture}(20,3.5)(-10,-0.5)
\put(0,2){\line(1,0){3}}
\put(0,1){\line(1,0){3}}
\put(0,1){\line(0,1){1}}
\put(1,1){\line(0,1){1}}
\put(2,1){\line(0,1){1}}
\put(3,1){\line(0,1){1}}
\put(3,0){\makebox(2,2){$ \subset$}}
\put(5,2){\line(1,0){4}}
\put(5,1){\line(1,0){4}}
\put(5,0){\line(1,0){1}}
\put(5,0){\line(0,1){2}}
\put(6,0){\line(0,1){2}}
\put(7,1){\line(0,1){1}}
\put(8,1){\line(0,1){1}}
\put(9,1){\line(0,1){1}}
\put(9,0){\makebox(2,2){$ \subset$}}
\put(11,2){\line(1,0){4}}
\put(11,1){\line(1,0){4}}
\put(11,0){\line(1,0){3}}
\put(11,0){\line(0,1){2}}
\put(12,0){\line(0,1){2}}
\put(13,0){\line(0,1){2}}
\put(14,0){\line(0,1){2}}
\put(15,1){\line(0,1){1}}
\put(16,0){\makebox(4,2){${\displaystyle {3 \brack 2}{3 \brack 1}}$.}}
\end{picture}
\\ By summing them up, we obtain
(\ref{qpol2}).
\end{example}
%%%%%%%%%%%%%%%%%%%%%%%%%%%%%%%%%%%%%%%%%%%%%%%%%%%%%%%%%%%%%%%%%%
\Proof (Proposition \ref{pro:ffkkp})
Let us note that the left hand side of the equation
(\ref{eq:ffkkp}) is nothing but the transition coefficient
{}from elementary symmetric function basis to Hall-Littlewood basis
$$
e_{\lambda}(x)=e_{\lambda_1}(x) \cdots e_{\lambda_{n}}(x)=
\sum_{\eta (\eta_1 \leq n)} \sum_{\mu}
K_{\eta' \lambda} K_{\eta \mu}(q) P_{\mu}(x;q),
$$
which follows from $e_{\lambda}(x)=\sum_{\eta} K_{\eta' \lambda} s_{\eta}(x)$
and $s_{\eta}(x)=\sum_{\mu} K_{\eta \mu}(q) P_{\mu}(x;q)$.
Then by using the formula (see, e.g.,\cite{Ma} eq.(3.2) p215)
\begin{eqnarray}
e_m(x) P_{\nu}(x;q)&=&\sum_{\mu} f^{\mu}_{\nu (1^m)} P_{\mu}(x;q),\cr
f^{\mu}_{\nu (1^m)}&=& \prod_{i \geq 1}
\left[ \begin{array}{c} \mu'_i -  \mu'_{i+1}
 \\   \mu'_i -  \nu'_i
\end{array} \right],
\label{eq:etop}
\end{eqnarray}
and therefore $f^{\mu}_{\nu,(1^m)}=0$ unless $\mu/\nu$ is a
vertical $m$-strip, we obtain
$$
e_{\lambda}(x)=\sum_{\{ \mu \}}
f^{\mu^{(n)}}_{\mu^{(n-1)} (1^{\lambda_n})}\cdots
f^{\mu^{(2)}}_{\mu^{(1)} (1^{\lambda_2})}
f^{\mu^{(1)}}_{\mu^{(0)} (1^{\lambda_1})} P_{\mu^{(n)}}(x;q),
$$
where the sum is taken over the partitions
$\{ \mu^{(0)}=\emptyset,\mu^{(1)}, \ldots, \mu^{(n)} \}.$
By putting ${\mu^{(a)}}'=\nu^{(a)}$ we have the desired
formula (\ref{eq:ffkkp}).
\qed

\medskip
\Proof (Proposition \ref{pro:ffkk})
The proof is essentially the same as the proof of
Proposition \ref{pro:ffkkp}.
The relevant transition coefficients are those
{}from complete symmetric function basis to Hall-Littlewood basis
$$
h_{\lambda}(x)=h_{\lambda_1}(x) \cdots h_{\lambda_{n}}(x)=
\sum_{\eta (l(\eta) \leq n)} \sum_{\mu}
K_{\eta \lambda} K_{\eta \mu}(q) P_{\mu}(x;q),
$$
which follows from $h_{\lambda}(x)=\sum_{\eta} K_{\eta \lambda} s_{\eta}(x)$
and $s_{\eta}(x)=\sum_{\mu} K_{\eta \mu}(q) P_{\mu}(x;q)$.
In this case, the formula corresponding to (\ref{eq:etop})
is the following,
\begin{eqnarray}
h_m(x) P_{\nu}(x;q)&=&\sum_{\mu} g^{\mu}_{\nu (m)} P_{\mu}(x;q),\cr
g^{\mu}_{\nu (m)}&=&q^{ \displaystyle{\sum_{i \geq 1}}
\left( \begin{array}{c} \mu'_i -  \nu'_i \\
2 \end{array} \right)}
\prod_{i \geq 1}
\left[ \begin{array}{c} \mu'_i -  \nu'_{i+1}
 \\   \nu'_i -  \nu'_{i+1}
\end{array} \right],
\label{eq:htop}
\end{eqnarray}
and $g^{\mu}_{\nu (m)}=0$ unless $\nu \subset \mu$
and $\vert \mu/\nu \vert=m$.
The equation (\ref{eq:htop}) can be proved as follows.

Because of the relation
$\sum_{l=0}^{m} (-1)^l e_{m-l}(x) h_{l}(x)=\delta_{m,0}$,
we have
\begin{equation}
\sum_{\mu}
f^{\mu}_{\nu (1^k)} g^{\lambda}_{\mu (l)} (-1)^l=\delta_{\lambda, \nu},
\label{eq:fgrel}
\end{equation}
where
$k=\vert \mu/\nu \vert$,
$l=\vert \lambda/\mu \vert$ and
the sum is taken for all $\mu$ such that
$\nu \subset \mu \subset \lambda$.
Since $\{g^{\mu}_{\nu (l)}\}$ are uniquely determined from
$\{f^{\mu}_{\nu (1^k)}\}$ by this relation(\ref{eq:fgrel}),
it is enough to show that the explicit formulae for
$g^{\mu}_{\nu (l)}$ (\ref{eq:htop}) and
$f^{\mu}_{\nu (1^k)}$ (\ref{eq:etop}) indeed satisfy this relation.
To this end, we
rewrite the left hand side of the relation (\ref{eq:fgrel})
as follows
\begin{equation}
\sum_{\nu\subset\mu\subset\lambda}(-1)^{\sum(\lambda_i'-\mu_i')}
q^{\sum\pmatrix{\lambda_i'-\mu_i'\cr 2}}\prod_{i\ge 1}
\left[\matrix{\lambda_i'-\mu_{i+1}'\cr \lambda_i'-\mu_i'}\right]
\left[\matrix{\mu_i'-\mu_{i+1}'\cr\mu_i'-\nu_i'}\right]
=\prod_{i \geq 1} \Phi_i(q),
\label{eq:phi}
\end{equation}
where
$$
\Phi_i(q)=\sum_{\nu_i'\le\mu_i'\le\lambda_i'}
(-1)^{\lambda_i'-\mu_i'}
q^{\pmatrix{\lambda_i'-\mu_i'\cr 2}}{(q)_{\lambda_{i-1}'-\mu_i'}\over
(q)_{\lambda_i'-\mu_i'}(q)_{\mu_i'-\nu_i'}(q)_{\nu_{i-1}'-\mu_i'}},
$$
and $(q)_{\lambda_0'-\mu_1'}/(q)_{\nu_0'-\mu_1'}=1$.
Precisely speaking, the sum in equation $\Phi_i$
should be taken over the $\mu'_i$ such as
$\max(\mu'_{i+1},\nu'_i) \leq \mu'_i \leq \lambda'_i$,
since $\mu$ is a partition.
However, if $\mu'_{i+1}>\nu'_i$, one can consider $\Phi_{i+1}$ vanishes
by understanding $1/(q)_{\nu'_i-\mu'_{i+1}}=0$,
hence one can simply write as above.

Consider at first $\Phi_1(q)$. We have
\begin{eqnarray}
(q)_{\lambda_1'-\nu_1'}\Phi_1(q)
&=&
\sum_{\nu'_1 \leq \mu'_1 \leq \lambda'_1}(-1)^{\lambda_1'-\mu_1'}
q^{\pmatrix{\lambda_1'-\mu_1'\cr 2}}
\left[\matrix{\lambda_1'-\nu_1'\cr \lambda_1'-\mu_1'}\right] \cr
&=&
\sum_{m=0}^{\lambda'_1-\nu'_1}
(-1)^m q^{\pmatrix{m\cr 2}}
\left[\matrix{\lambda_1'-\nu_1'\cr m}\right]
=\delta_{\lambda_1',\nu_1'}.
\end{eqnarray}
The last equality follows from the $q$--binomial theorem
$$\sum_{m=0}^N(-z)^mq^{\pmatrix{m\cr 2}}\left[\matrix{N\cr m}\right]=
\prod_{i=1}^N(1-q^{i-1}z).
$$
Thus we have $\Phi_1=\delta_{\lambda'_1, \nu'_1}$.
Similarly, under the condition $\lambda_1'=\nu_1'$, we obtain
$\Phi_2=\delta_{\lambda'_2, \nu'_2}$, since
$$(q)_{\lambda_2'-\nu_2'}\Phi_2(q)=
\sum_{m=0}^{\lambda'_2-\nu'_2}(-1)^m
q^{\pmatrix{m\cr 2}}\left[\matrix{\lambda_2'-\nu_2'\cr m}\right]=
\delta_{\lambda_2',\nu_2'}.
$$
Repeating these arguments we see that the product (\ref{eq:phi})
is equal to $\delta_{\lambda,\nu}$.
This proves (\ref{eq:fgrel}) and
hence the Proposition \ref{pro:ffkk}.
\qed

Essentially the same proof is also available in \cite{Ki3}
as well as discussions on several other aspects.
%\medskip
%%%%%%%%%%%%%%%%%%%%%%%%%%%%%%%%%%%%%%%%%%%%%%%%%%%%%%%%%%%%%%%%%%%

The classically restricted 1dsum
${X}_\mu(\la)$
(\ref{xbarbyk}) is nothing but the Kostka-Foulkes polynomial.
Its fermionic formula reads as
\begin{proposition}[{\rm \cite{KR}}]\label{pro:ffk}
For any partitions $\la$ and $\mu$ such that
$\vert \la \vert = \vert \mu \vert$ and $l(\la) \le n$
we have
\begin{eqnarray}
K_{\la \mu}(q) & = & \sum_{\{m \}} q^{c(\{ m \})}
\prod_{{\scriptstyle 1 \le a \le n-1} \atop
   {\scriptstyle i \ge 1}}
\left[ \begin{array}{c} p^{(a)}_i +  m^{(a)}_i
 \\   m^{(a)}_i \end{array} \right],\label{eq:ffk1}\\
c(\{ m  \}) & = & n(\mu) + \frac{1}{2}\sum_{1 \le a, b \le n-1} C^{(n)}_{a b}
\sum_{j, k \ge 1} \mbox{min}(j, k) m^{(a)}_j m^{(b)}_k \nonumber\\
&& \qquad\quad - \sum_{j, k \ge 1} \mbox{min}(j, \mu_k)
m^{(1)}_j,\label{eq:ffk2}\\
p^{(a)}_i & = & \delta_{a 1} \sum_{k \ge 1} \mbox{min}(i, \mu_k)
- \sum _{b=1}^{n-1} C^{(n)}_{a b} \sum_{k \ge 1} \mbox{min}(i,k) m^{(b)}_k,
\label{eq:ffk3}
\end{eqnarray}
where the sum $\sum_{\{ m \}}$ is taken over
$\{ m^{(a)}_i \in {\bf Z}_{\ge 0} \mid 1 \le a \le n-1, \, i \ge 1 \}$
satisfying $p^{(a)}_i \ge 0$ for $1 \le a \le n-1, i \ge 1$, and
\begin{equation}\label{eq:ffk4}
\sum_{i \ge 1} i\, m^{(a)}_i = \la_{a+1} + \la_{a+2} + \cdots + \la_n
\qquad \mbox{ for } \ 1 \le a \le n-1.
\end{equation}
\end{proposition}
This formula is originated in the Bethe ansatz.
When $q=1$ it counts the
multiplicity of the $\la$-representation in the tensor product
of the $\mu_i$-fold symmetric tensor representations
$( 1 \le i \le l(\mu))$.
Kostka-Foulkes polynomials
$K_{\xi \eta}(q)$ can also be realized through the dual picture
corresponding to the 1dsum $X'_\eta(\xi')$, namely,
the ($q$-)multiplicity of the $\xi'$-representation in the tensor product
of the $\eta_i$-fold antisymmetric tensor representations
$( 1 \le i \le l(\eta))$.
Such a duality will be discussed in a more general setting
in \cite{KS}.
{}From the results therein we have
another fermionic formula as
\begin{proposition}[{\rm \cite{KS}}]\label{pro:affk}
For any partitions $\xi$ and
$\eta = ((n-1)^{\zeta_{n-1}} \cdots 1^{\zeta_1})$ such that
$\vert \xi \vert = \vert \eta \vert$ and $\xi_1 \le n$, we have
\begin{eqnarray}
K_{\xi \eta}(q) & = & \sum_{\{\hat{m} \}} q^{\hat{c}(\{ \hat{m} \})}
\prod_{{\scriptstyle 1 \le a \le n-1} \atop
   {\scriptstyle i \ge 1}}
\left[ \begin{array}{c} \hat{p}^{(a)}_i +  \hat{m}^{(a)}_i
 \\   \hat{m}^{(a)}_i \end{array} \right],\label{eq:affk1}\\
\hat{c}(\{ \hat{m}  \}) & = & \frac{1}{2}
\sum_{1 \le a, b \le n-1} C^{(n)}_{a b}
\sum_{j, k \ge 1} \mbox{min}(j, k)
\hat{m}^{(a)}_j \hat{m}^{(b)}_k,\label{eq:affk2}\\
\hat{p}^{(a)}_i & = & \zeta_a
- \sum _{b=1}^{n-1} C^{(n)}_{a b} \sum_{k \ge 1}
\mbox{min}(i,k) \hat{m}^{(b)}_k,
\label{eq:affk3}
\end{eqnarray}
where the sum $\sum_{\{ \hat{m} \}}$ is taken over
$\{ \hat{m}^{(a)}_i \in {\bf Z}_{\ge 0} \mid 1 \le a \le n-1, \, i \ge 1 \}$
satisfying $\hat{p}^{(a)}_i \ge 0$ for
$1 \le a \le n-1, i \ge 1$, and
\begin{equation}\label{eq:affk4}
\sum_{i \ge 1} i\,  \hat{m}^{(a)}_i = \sum_{b=1}^{n-1}
\mbox{min}(a, b) \zeta_b - (\xi'_1 + \cdots + \xi'_a)
\qquad \mbox{ for } \ 1 \le a \le n-1.
\end{equation}
\end{proposition}

As for the restricted 1dsums
$X^{(l)}_\mu(\la)$ and $X^{(l)'}_\mu(\la)$,
fermionic formulae are yet conjectural in general and
not yet available for arbitrary $\la$.
($\goth{sl}_{\, 2}$ case is the exception,
see \cite{BMS}, \cite{Ki} and references therein.)
Here we shall only deal with the vacuum case
$\la = l\Lambda_0$ corresponding to the partitions
$((\frac{\vert \mu \vert}{n})^n)$ and
$(n^{\frac{\vert \mu \vert}{n}})$, respectively.
See Section \ref{discussionhata} for more general
$\la$ cases.
To present our conjecture we prepare two expressions
$F^{(l)}_\mu(q)$ and $F^{(l)'}_\eta(q)$.
The first one is defined for partitions
$\mu$ satisfying $\vert \mu \vert \equiv 0$ mod $n$,
$\mu_1 \le l$ and reads
\begin{eqnarray}
F^{(l)}_{\mu}(q) & = & \sum_{\{m \}} q^{c_l(\{ m \})}
\prod_{{\scriptstyle 1 \le a \le n-1} \atop
   {\scriptstyle 1 \le i \le l-1}}
\left[ \begin{array}{c} p^{(a)}_i +  m^{(a)}_i
 \\   m^{(a)}_i \end{array} \right],\label{eq:ffrk1}\\
c_l(\{ m  \}) & = & n(\mu) +
\frac{1}{2}\sum_{1 \le a, b \le n-1} C^{(n)}_{a b}
\sum_{1 \le j, k \le l} \mbox{min}(j, k) m^{(a)}_j m^{(b)}_k \nonumber\\
&& \qquad\quad - \sum_{k \ge 1} \sum_{j=1}^l \mbox{min}(j, \mu_k)
m^{(1)}_j,\label{eq:ffrk2}\\
p^{(a)}_i & = & \delta_{a 1} \sum_{k \ge 1} \mbox{min}(i, \mu_k)
- \sum _{b=1}^{n-1} C^{(n)}_{a b}
\sum_{k=1}^l \mbox{min}(i,k) m^{(b)}_k,
\label{eq:ffrk3}
\end{eqnarray}
where the sum $\sum_{\{ m \}}$ is taken over
$\{ m^{(a)}_i \in {\bf Z}_{\ge 0} \mid 1 \le a \le n-1,
\, 1 \le i \le l \}$
satisfying $p^{(a)}_i \ge 0$ for
$1 \le a \le n-1,\, 1 \le i \le l-1$, and
\begin{equation}\label{eq:ffrk4}
\sum_{i=1}^l i\,  m^{(a)}_i = \frac{n-a}{n}\vert \mu \vert
\qquad \mbox{ for } \ 1 \le a \le n-1.
\end{equation}
This is very similar to the fermionic form appearing in
Proposition \ref{pro:ffk}
with $\la = ((\frac{\vert \mu \vert}{n})^n)$.
In fact (\ref{eq:ffrk1})-(\ref{eq:ffrk4}) correspond to truncating
the indices of $m^{(a)}_i$ and $p^{(a)}_i$ in
(\ref{eq:ffk1})-(\ref{eq:ffk4})
to the range $1 \le i \le l$.
The second one $F^{(l)'}_\eta(q)$ is a similar analogue of
the fermionic form in Proposition \ref{pro:affk}.
It is defined for the partitions of the form
$\eta = ((n-1)^{\zeta_{n-1}} \cdots 1^{\zeta_1})$
such that $\vert \eta \vert \equiv 0$ mod $n$.
\begin{eqnarray}
F^{(l)'}_\eta(q) & = &
\sum_{\{\tilde{m} \}} q^{\tilde{c}_l(\{ \tilde{m} \})}
\prod_{{\scriptstyle 1 \le a \le n-1} \atop
   {\scriptstyle 1 \le i \le l-1}}
\left[ \begin{array}{c} \tilde{p}^{(a)}_i +  \tilde{m}^{(a)}_i
 \\   \tilde{m}^{(a)}_i \end{array} \right],\label{eq:affrk1}\\
\tilde{c}_l(\{ \tilde{m}  \}) & = & \frac{1}{2}
\sum_{1 \le a, b \le n-1} C^{(n)}_{a b}
\sum_{1 \le j, k \le l} \mbox{min}(j, k)
\tilde{m}^{(a)}_j \tilde{m}^{(b)}_k,\label{eq:affrk2}\\
\tilde{p}^{(a)}_i & = & C^{(l)^{-1}}_{1 i}\zeta_a
- \sum _{b=1}^{n-1} C^{(n)}_{a b} \sum_{k=1}^{l-1}
C^{(l)^{-1}}_{i k} \tilde{m}^{(b)}_k,
\label{eq:affrk3}
\end{eqnarray}
where the sum $\sum_{\{ \tilde{m} \}}$ is taken over
$\{ \tilde{m}^{(a)}_i \in {\bf Z}_{\ge 0}
\mid 1 \le a \le n-1, \, 1 \le i \le l \}$
satisfying $\tilde{p}^{(a)}_i \ge 0$ for
$1 \le a \le n-1, \, 1 \le i \le l-1$, and
\begin{equation}\label{eq:affrk4}
\sum_{i=1}^l  i\,  \tilde{m}^{(a)}_i = \sum_{b=1}^{n-1}
C^{(n)^{-1}}_{a b} \zeta_b
\qquad \mbox{ for } \ 1 \le a \le n-1.
\end{equation}
This is a truncation of (\ref{eq:affk1})--(\ref{eq:affk4})
with $\xi = (n^{\frac{\vert \eta \vert}{n}})$.
Note that the variables $\tilde{m}^{(a)}_l \, (1 \le a \le n-1)$
are specified and appear only in (\ref{eq:affrk4}) to impose
the condition $\tilde{m}^{(a)}_l \in {\bf Z}_{\ge 0}$.
Now our conjecture is stated as
\begin{conjecture}\label{con:ffrk}
For a partition $\mu$ such that
$\vert \mu \vert \equiv 0 \mbox{ mod } n$ and $\mu_1 \le l$ we have
\begin{equation}\label{con:ffrk1}
X^{(l)}_{\mu}(l\La_0)  =  F^{(l)}_\mu(q).
\end{equation}
For a partition $\eta$ such that
$\eta_1 \le n-1$ and
$\vert \eta \vert \equiv 0 \mbox{ mod } n$ we have
\begin{equation}\label{con:ffrk2}
X^{(l)'}_{\eta}(l\La_0)  =  F^{(l)'}_\eta(q).
\end{equation}
\end{conjecture}
This is also based on the Bethe ansatz for RSOS models,
in which an analogous truncation
by the level takes place.
In Section \ref{discussionhata} we will present a more general conjecture
on $X^{(l)'}_\eta(\la)$ for $\la$ non vacuum type.
%
%
%******************subsection: FF of the limits **************
%
\subsection{Fermionic formulae of the limits}\label{sec:fflimit}
By taking the large $\mu$ limit of the 1dsums in Section \ref{ssec4.1},
we obtain explicit $q$-series formulae of the characters in
Proposition \ref{1dsumlimit}.
For the unrestricted 1dsum $g_\mu(\la)$ we consider two
particular limits.
The first one is $\mu = (l^L) \cup \nu$, where
$L \rightarrow \infty$ under the condition
$L \equiv r\ (0 \le r \le n-1)$ mod $n$ and $\nu$ is a fixed
finite partition with $\nu_1 < l$.
By (\ref{vcal}) this corresponds to the module
${\cal V} = \oplus_{p \in {\cal H}(l\Lambda_r, \nu)}
V(l\Lambda_r + af(\wt p) -
(E(p) - \overline{E}(l\La_r, \nu))\delta)$.
See Example \ref{ex:calB}.
By calculating the above limit of (\ref{eq:ffkk}) one can show
\begin{proposition}\label{pro:kklim1}
Let $0 \le r \le n-1$ and
${\cal V} = \oplus_{p \in {\cal H}(l\Lambda_r, \nu)}
V(l\Lambda_r + af(\wt p) -
(E(p) - \overline{E}(l\La_r, \nu))\delta)$.
For $\la = (\la_1-\la_2)\overline{\Lambda}_1 + \cdots
+ (\la_{n-1} - \la_n)\overline{\Lambda}_{n-1} \in
\overline{P}$, we have
\begin{eqnarray}\label{eq:kklim1}
q^\Delta c^{\cal V}_\la(q)
%\sum_i(\mbox{dim}{\cal V}_{\la + l\La_0-i\delta}) q^i
& = &
\frac{1}{(q)_{\infty}^{n-1}} \sum_{\{ m\} }
\frac{q^{C(\{ m \} )}}
{\prod_{{\scriptstyle 1 \le a \le n-1} \atop
   {\scriptstyle 1 \le i \le l-1}}
(q)_{m_i^{(a)}}}, \\
C(\{ m \} ) &=&
\frac{1}{2} \sum_{a,b=1}^{n-1} \sum_{j,k=1}^{l-1}
C^{(n)}_{ab} C^{(l)-1}_{jk} m^{(a)}_j m^{(b)}_k \nonumber \\
& & - \sum_{j,k=1}^{l-1}C^{(l)^{-1}}_{j k}
(\nu'_j - \nu'_{j+1}) m^{(n-1)}_{k},
\end{eqnarray}
where
$\Delta = \frac{1}{2}\sum_{j=1}^l(\La_{r+\nu'_j} \vert \La_{r+\nu'_j})
- \frac{(\lambda \vert \lambda)}{2l}
-\frac{n-1}{2n}\sum_{j,k=1}^{l-1}
C^{(l)^{-1}}_{j k}(\nu'_j - \nu'_{j+1})(\nu'_k - \nu'_{k+1})$.
The sum $\sum_{\{ m \}}$ is taken over
$\{ m^{(a)}_i \in {\bf Z}_{\ge 0} \mid 1 \le a \le n-1, \,
1 \le i \le l-1 \}$ satisfying
\begin{equation}
\sum_{i=1}^{l-1} i\, m^{(a)}_i \equiv
\sum_{b=1}^a\left(\la_{b} -
\frac{\vert \la \vert - \vert \nu \vert - lr}{n} \right)
\ \  \mbox{ mod }\mbox{ }  l
\qquad \mbox{ for } \ 1 \le a \le n-1.
\end{equation}
\end{proposition}

\Proof
Let us rewrite (\ref{eq:ffkk})-(\ref{eq:ffkkphi})
in terms of the
variables $m^{(a)}_i = \nu^{(a)}_i - \nu^{(a)}_{i+1},
\, 1 \le i \le l = \mu_1$:
\begin{eqnarray}
\sum_{\eta \, (l(\eta) \le n)} K_{\eta \la}K_{\eta \mu}(q) & = &
\sum_{\{m \}} q^{\psi(\{m \})}
\prod_{{\scriptstyle 1 \le a \le n-1} \atop
   {\scriptstyle 1 \le i \le l}}
\left[ \begin{array}{c} p^{(a)}_i + m^{(a)}_{i}
 \\   m^{(a)}_i \end{array} \right], \label{proofkk} \\
p^{(a)}_i & = & \sum_{j=i}^l(m^{(a+1)}_j - m^{(a)}_j),\label{proofkkp}\\
\psi(\{m \}) & = &  n(\mu) +  \frac{(\la \vert \la)}{2l} -
\frac{n-1}{2nl}\vert \mu \vert^2 \nonumber \\
& &  + \frac{1}{2} \sum_{a,b=1}^{n-1} \sum_{j,k=1}^{l-1}
C^{(n)}_{ab} C^{(l)-1}_{jk} m^{(a)}_j m^{(b)}_k \nonumber \\
& & - \sum_{j,k=1}^{l-1}C^{(l)^{-1}}_{j k}
(\mu'_j - \mu'_{j+1}) m^{(n-1)}_{k},\label{proofkkpsi}
\end{eqnarray}
where $m^{(n)}_i = \mu'_i - \mu'_{i+1}$.
In view of (\ref{eq:ffkknu}) the sum $\sum_{\{ m \}}$ in
(\ref{proofkk}) runs over
$\{ m^{(a)}_i \in {\bf Z}_{\ge 0} \mid 1 \le a \le n-1, 1 \le i \le l \}$
obeying the condition
\begin{equation}\label{proofkkwa}
\sum_{i=1}^l i\, m^{(a)}_i = \la_1 + \cdots + \la_a.
\end{equation}
In deriving (\ref{proofkkpsi}) we have eliminated
$m^{(a)}_l$ for $1 \le a \le n-1$ by (\ref{proofkkwa})
and used $(\la \vert \la) = \sum_{a=1}^n \la_a^2 - \vert \la \vert^2/n$.
The limit $L \rightarrow \infty$ of (\ref{proofkk}) is to be expanded from
such $m^{(a)}_i$ that attains the minimum of
$\psi$ (\ref{proofkkpsi}).
It occurs around
$m^{(a)}_i = \frac{a}{n}(L-r)\delta_{i l}$.
Thus $\forall p^{(a)}_i \rightarrow \infty$
in the limit $L \rightarrow \infty$ hence
the product of the $q$-binomial coefficients tends to
$\left( (q)^{n-1}_\infty \prod_{a=1}^{n-1}\prod_{i=1}^{l-1}
(q)_{m^{(a)}_i}\right)^{-1}$.
The powers of $q$ can be adjusted by noting
$\mu'_j - \mu'_{j+1} = \nu'_j - \nu'_{j+1}$ for the shape
$\mu = (l^L) \cup \nu$ in (\ref{proofkkpsi}) and
verifying
$-\overline{E}(l\La_0, \mu) + \Delta +
\psi(\{m \})\vert_{m^{(a)}_i \rightarrow
m^{(a)}_i + \frac{a}{n}(L-r)\delta_{i l}} = C(\{m\})$ by
using the ground state energy (\ref{gsenergy}). \qed

For $r=0$ and $\nu = \emptyset$ one has
${\cal V} = V(l\La_0)$, in which case
the above formula was conjectured in \cite{KNS},
announced in \cite{FS} and proved in \cite{Ge2}.
When $r=n-1$ and $\nu = (s)$ it agrees with $j=n$ case of eq.(5.7)
in \cite{Ge2}.
However the above result does not cover the $V((l-s)\La_0 + s\La_j)$
case for general $j$ obtained in \cite{Ge2}.
\begin{remark}\label{rem:tba}
Although the fermionic form (\ref{proofkk}) looks
similar to that in Proposition \ref{pro:ffk},
its interpretation in terms of the Bethe ansatz is yet
unknown to us.
\end{remark}

For an $U_q(\widehat{\goth{sl}}_{\, n})$-module $M$,
let ${\sl ch}\, M$ denote its character:
$$
{\sl ch}\, M = \sum_\la (\mbox{dim} M_\la)\, e^\la,
$$
where $q = e^{-\delta}$ and the sum runs over all the affine weights
$\la \in P$.
Then we make
\begin{remark}\label{rem:fsg}
Substituting $r = 0, \nu = \emptyset$ case of
Proposition \ref{pro:kklim1} into eq.(12.7.12) in \cite{Ka},
one can express the character of the
level $l$ vacuum module ${\cal V} = V(l\La_0)$ as
\begin{eqnarray}\label{eq:fsg}
e^{-l\La_0} {\sl ch } {V(l\La_0)} & = &
\frac{1}{(q)_{\infty}^{n-1}} \sum_{\{ m\} }
\frac{q^{\overline{C}_l(\{ m \} )}}
{\prod_{{\scriptstyle 1 \le a \le n-1} \atop
   {\scriptstyle 1 \le i \le l-1}}
(q)_{m_i^{(a)}}} \ e^{\beta_l(\{m\})}, \\
\overline{C}_l(\{ m \} ) & = &
\frac{1}{2} \sum_{a,b=1}^{n-1} \sum_{j,k=1}^{l}
C^{(n)}_{ab} \mbox{min}(j,k) m^{(a)}_j m^{(b)}_k, \\
\beta_l(\{m \}) & = & \sum_{a=1}^{n-1}\sum_{i=1}^l i\, m^{(a)}_i \alpha_a.
\end{eqnarray}
Here the sum $\sum_{\{ m \} }$ runs over
$\{m^{(a)}_i \in {\bf Z}_{\ge 0} \mid 1 \le a \le n-1,
1 \le i \le l-1 \}$ and
$\{ m^{(a)}_l \in {\bf Z} \mid 1 \le a \le n-1 \}$.
This agrees with a formula announced in \cite{FS}
and proved in \cite{Ge1}.
\end{remark}
%
%********** prop: KK limit type 2 *****************
%
The second limit we consider of $g_\mu(\la)$ given by (\ref{gbyk}) and
Proposition \ref{pro:ffkk} is
$\mu = (l_1^{L_1}) + \cdots + (l_s^{L_s})$ in the notation of
Section \ref{subsec:limit}, in which
$L_J - L_{J+1} \rightarrow \infty$ for $1 \le J \le s$
keeping $L_J \equiv r_J$ mod $n$ fixed
($0 \le r_J \le n-1, \, L_{s+1} = r_{s+1} = 0$).
A similar calculation to Proposition \ref{pro:kklim1}
leads to
\begin{proposition}\label{pro:kklim2}
Fix $l_1, \ldots,l_s \in {\bf Z}_{\ge 0}$ and
$0 \le r_1, \ldots, r_s  \le n-1$.
Set $\vert r \vert = r_1 + \cdots + r_s, \, l= l_1 + \cdots + l_s$ and
$M_J = l_1 + l_2 + \cdots + l_J$.
For $\la = (\la_1-\la_2)\overline{\Lambda}_1 + \cdots
+ (\la_{n-1} - \la_n)\overline{\Lambda}_{n-1} \in
\overline{P}$ and the tensor product module
${\cal V} = \otimes_{J=1}^s V(l_J \Lambda_{r_J})$,
we have
\begin{eqnarray}\label{eq:kklim2}
q^\Delta c^{\cal V}_\la(q)
%\sum_i(\mbox{dim}{\cal V}_{\la + l\La_0- i \delta}) q^i
& = &
\frac{1}{(q)_{\infty}^{(n-1)s}} \sum_{\{ m\} }
\frac{q^{C(\{ m \} )}}{\prod_{a=1}^{n-1}
\prod_{{\scriptstyle 1 \le i \le l-1} \atop
   {\scriptstyle i \notin \{M_1, \ldots, M_{s-1}\}}}
(q)_{m_{i}^{(a)}}} \\
C(\{ m \} ) &=&
\frac{1}{2} \sum_{a,b=1}^{n-1} \sum_{j,k=1}^{l-1}
C^{(n)}_{ab} C^{(l)^{-1}}_{jk} m^{(a)}_j m^{(b)}_k \nonumber \\
&& -\sum_{j=1}^{l-1}\sum_{J=1}^{s-1}
C^{(l)^{-1}}_{j\, M_J} (r_J-r_{J+1})m^{(n-1)}_j,
\end{eqnarray}
where $\Delta =
\sum_{J=1}^s\frac{l_J}{2}(\Lambda_{r_J}\vert \Lambda_{r_J}) -
\frac{(\la \vert \la)}{2l} -
\frac{n-1}{2n}\sum_{J,K=1}^{s-1}
C^{(l)^{-1}}_{M_J\, M_K}(r_J-r_{J+1})(r_K-r_{K+1})$.
The sum $\sum_{\{ m \}}$ runs over
$\{ m^{(a)}_i \mid 1 \le a \le n-1, \,
1 \le i \le l-1 \}$ satisfying
\begin{eqnarray}
m_i^{(a)} & \in & \left\{
\begin{array}{ll}
{\bf Z} & \quad i \in \{M_1, \ldots, M_{s-1}\}\\
{\bf Z}_{\ge 0} & \quad \mbox{ otherwise}
\end{array} \right. ,\label{hanni}\\
\sum_{i=1}^{l-1} i\, m^{(a)}_i & \equiv &
\sum_{b=1}^a\left(\la_b -
\frac{| \la |-\sum_{J=1}^{s}l_J r_J}{n}\right)
\  \mbox{mod }\, l \ \mbox{for } 1 \le a \le n-1.
\label{modcond}
\end{eqnarray}
\end{proposition}
\Proof
Again we start with the expression
(\ref{proofkk})--(\ref{proofkkpsi}) with
$\mu = (l_1^{L_1}) + \cdots + (l_s^{L_s})$.
This time the minimum of $\psi$ (\ref{proofkkpsi})
is attained around
$m^{(a)}_{i, 0} = \frac{a}{n}(L_J - L_{J+1} - r_J + r_{J+1})$
for $i = M_J \, (1 \le J \le s)$,
$ = 0$ for $i \notin \{M_1, \ldots, M_s \}$.
In the limit $L_J - L_{J+1} \rightarrow \infty$,
$m^{(a)}_{i,0}$ tends to infinity for
$i \in \{M_1, \ldots, M_s \}$ and so does $\forall p^{(a)}_i$.
Thus the product of the $q$-binomial coefficients converges to
$\left((q)_{\infty}^{(n-1)s}
\prod_{a=1}^{n-1}
\prod_{i \in \{1, \ldots, l-1 \} \setminus \{M_1, \ldots, M_{s-1}\}}
(q)_{m_{i}^{(a)}}\right)^{-1}$.
After the replacement
$m^{(a)}_i \rightarrow m^{(a)}_i + m^{(a)}_{i, 0}$,
the new variables $m^{(a)}_i$ are to satisfy
(\ref{hanni})--(\ref{modcond}).
It remains to check
$-\overline{E}(l\La_0, \mu) + \Delta +
\psi(\{m \})\vert_{m^{(a)}_i \rightarrow m^{(a)}_i + m^{(a)}_{i, 0}}
= C(\{m\})$ for $\mu = (l_1^{L_1}) + \cdots + (l_s^{L_s})$.
This is straightforward by using the explicit forms
(\ref{gsenergy}) and (\ref{proofkkpsi}).
\qed

By the definition the module ${\cal V}$ remains unchanged under any
permutations of $(l_J, r_J), \, (J = 1 , \ldots, s)$.
Note however that such symmetry is not manifest in the RHS
of (\ref{eq:kklim2})--(\ref{hanni}).
When $r_1 = \cdots = r_s = 0$, Proposition \ref{pro:kklim2}
can also be shown by decomposing
the product of $e^{-l_J\La_0}{\sl ch } {V(l_J\La_0)}$ over $1 \le J \le s$
given in Remark \ref{rem:fsg}.
This will be explained more precisely in
Section \ref{discussionkuni} together with a
conjectural extension of Proposition \ref{pro:kklim2}
to an arbitrary non twisted affine Lie algebra $X^{(1)}_n$.

Let us proceed to the limit of the classically restricted 1dsum
$X_\mu(\la)$.
Here we shall exclusively consider the situation
$\mu = (l^L)$, in which $L \rightarrow \infty$
under the condition $L \equiv 0$ mod $n$.
{}From Example \ref{ex:calB} (1) this is related to the vacuum module
${\cal V} = V(l\La_0)$.
Our task is to compute the limit of the Kostka-Foulkes polynomial
given in Proposition \ref{pro:ffk}.
Combining the result with (\ref{xbarbyk}) and
(\ref{xbarlimit})  we obtain
%
%********** prop: K limit  *****************
%
\begin{proposition}\label{pro:klim}
For any $l \in {\bf Z}_{\ge 1}$ and $\la \in \overline{P}^+$
such that $\vert \la \vert \equiv 0$ mod $n$, we have
\begin{eqnarray}
b^{V(l\La_0)}_\la(q)
%\sum_i [V(l\La_0): \la + l\La_0 - i \delta]_{cl}\, q^i
& = &
\sum_\eta
\frac{K_{\xi \eta}(q)\,
F^{(l)'}_\eta(q)}
{(q)_{\zeta_1}\cdots (q)_{\zeta_{n-1}} },\label{eq:klim1}\\
\xi  =  (n^{\frac{\vert \eta \vert - \vert \la \vert}{n}})\cup \la',
\qquad
\eta & = & \left( (n-1)^{\zeta_{n-1}},\ldots ,1^{\zeta_{1}} \right),
\label{eq:klim2}\\
\nonumber
\end{eqnarray}
where the sum $\sum_\eta$ runs over the partitions
$\eta$ satisfying $\eta_1 \le n-1$ as above and
$|\eta| \equiv 0$ mod $n$.
\end{proposition}
\Proof
We start with the expression
(\ref{eq:ffk1})--(\ref{eq:ffk4})
with $\mu = (l^L)$.
The limit is to be expanded around the minimum of
$c(\{m\})$ (\ref{eq:ffk2}).
This takes place at
$m^{(a)}_{i,0} = \frac{n-a}{n}L\delta_{i l}$, which is tending to
infinity.
Thus under the identification $\zeta_a = p^{(a)}_l$ the factor
$\prod_{a=1}^{n-1} \left[ \begin{array}{c} p^{(a)}_l +  m^{(a)}_l
 \\   m^{(a)}_l \end{array} \right]$
in (\ref{eq:ffk1}) gives rise to
$\left((q)_{\zeta_1} \cdots (q)_{\zeta_{n-1}}\right)^{-1}$
as $L \rightarrow \infty$.
After the shift
$m^{(a)}_i \rightarrow m^{(a)}_i + m^{(a)}_{i, 0}$, the relations
(\ref{eq:ffk4}), (\ref{eq:ffk3}) and its $i=l$ case become
\begin{eqnarray}
\sum_{k \ge 1} k\, m^{(a)}_k & = &
\frac{a L l}{n} - \la_1 - \cdots - \la_a,\label{msum}\\
p^{(a)}_i & = & - \sum_{b=1}^{n-1} C^{(n)}_{a b}
\sum_{k \ge 1} \mbox{min}(i,k)m^{(b)}_k, \label{shiftedp}\\
\sum_{b=1}^{n-1}C^{(n)^{-1}}_{a b} \zeta_b & = &
- \sum_{k \ge 1} \mbox{min}(l,k) m^{(a)}_k. \label{zetadef}
\end{eqnarray}
Eliminating $m^{(a)}_l$ with (\ref{zetadef}),
one rewrites (\ref{shiftedp}) as
$$
p^{(a)}_i =
\left\{
\begin{array}{ll}
\zeta_a - \sum_{b=1}^{n-1}C^{(n)}_{a b} \sum_{k > l}
\mbox{min}(i-l,k-l)m^{(b)}_k & \quad i > l\\
C^{(l)^{-1}}_{1 \ l-i}\zeta_a -
\sum_{b=1}^{n-1}C^{(n)}_{a b}\sum_{k=1}^{l-1}
C^{(l)^{-1}}_{i k}m^{(b)}_k& \quad 1 \le i < l
\end{array} \right. .
$$
By setting
$p^{(a)}_i = \hat{p}^{(a)}_{i-l}$ and
$m^{(a)}_i = \hat{m}^{(a)}_{i-l}$, the $i > l$ case
in the above coincides with (\ref{eq:affk3}).
So does the $1 \le i < l$ case with (\ref{eq:affrk3})
under $p^{(a)}_i = \tilde{p}^{(a)}_{l-i}$ and
$m^{(a)}_i = \tilde{m}^{(a)}_{l-i}$.
Thus we are going to extract
$K_{\xi \eta}(q)$ in Proposition \ref{pro:affk} from
``$\prod_{a=1}^{n-1} \prod_{i>l}$ part''  and
$F^{(l)'}_\eta(q)$ in (\ref{eq:affrk1})--(\ref{eq:affrk4}) from
``$\prod_{a=1}^{n-1} \prod_{1 \le i < l}$ part''.
Actually the equality
$-\overline{E}(l\La_0,(l^L)) + c(\{m^{(a)}_i + m^{(a)}_{i,0}\})
= \hat{c}(\{\hat{m}\}) + \tilde{c}_l(\{\tilde{m}\})$
is valid among the quadratic forms
(\ref{eq:ffk2}), (\ref{eq:affk2}) and (\ref{eq:affrk2}).
It remains to check
(i) (\ref{eq:affk4}) for $\xi$ given by (\ref{eq:klim2}),
(ii) $\tilde{m}^{(a)}_l$ defined by (\ref{eq:affrk4}) is
a non-negative integer,
(iii) $\eta$ defined by (\ref{eq:klim2}) and (\ref{zetadef})
satisfies $\vert \eta \vert \equiv 0$ mod $n$.
To show (i),  combine (\ref{zetadef}) and (\ref{msum}) to see
$\sum_{k > l}(k-l)m^{(a)}_k = \sum_{i \ge 1} i \hat{m}^{(a)}_i
= \sum_{b=1}^{n-1}C^{(n)^{-1}}_{a b}\zeta_b +
\frac{a l L}{n} - \la_1 - \cdots - \la_a$.
Since $\xi'_a = \frac{\vert \eta \vert - Ll}{n} + \la_a$
{}from (\ref{eq:klim2}), the last quantity coincides with the
RHS of (\ref{eq:affk4}).
To show (ii), use (\ref{eq:affrk4}), (\ref{zetadef}) and
$\tilde{m}^{(a)}_{l-i} = m^{(a)}_i$ to see
$\tilde{m}^{(a)}_l = \frac{1}{l}\left(
\sum_{b=1}^{n-1}C^{(n)^{-1}}_{a b} \zeta_b -
\sum_{i=1}^{l-1}(l-i)m^{(a)}_i\right) =
- \sum_{i \ge 1} m^{(a)}_i \in {\bf Z}$.
Moreover from $i=1$ case of (\ref{shiftedp}), one has
$- \sum_{i \ge 1} m^{(a)}_i
= \sum_{b=1}^{n-1}C^{(n)^{-1}}_{a b} p^{(b)}_1 \ge 0$
because of $\forall C^{(n)^{-1}}_{a b} > 0$ and
$\forall p^{(b)}_1 \ge 0$.
To show (iii), use (\ref{zetadef}) to see
$\vert \eta \vert = \sum_{a=1}^{n-1} a \zeta_a =
-n \sum_{k \ge 1}\mbox{min}(l,k)m^{(n-1)}_k \equiv 0$
mod $n$.
\qed

By admitting (\ref{con:ffrk2}) Proposition \ref{pro:klim} confirms
the vacuum module case $\La = l\La_0$ of
\begin{conjecture}[{\rm spinon character formula, \cite{NY2}}]\label{spinonch1}
For any $\La \in (P^+)_l$ and $\la \in \overline{P}^+$
such that $\La \equiv \la + l\La_0$ mod $Q$, we have
\begin{equation}
b^{V(\La)}_\la(q)
%\sum_i [V(\La): \la + l\La_0 - i \delta]_{cl}\, q^i
=
\sum_\eta
\frac{X'_{\eta}(\la)\,
X^{(l)'}_\eta(\La)}
{(q)_{\zeta_1}\cdots (q)_{\zeta_{n-1}} },
\end{equation}
where the sum $\sum_\eta$ runs over the partitions of the form
$\eta = \left( (n-1)^{\zeta_{n-1}},\ldots ,1^{\zeta_{1}} \right)$
satisfying $\vert \eta \vert \equiv \vert \la \vert$ mod $n$.
\end{conjecture}

Finally we turn to a limit related to the restricted 1dsum
$X^{(l)}_\mu(l\La_0)$.
By a parallel calculation with Proposition \ref{pro:klim}
one can derive
%
%********** prop: Restricted K limit  *****************
%
\begin{proposition}\label{pro:rklim}
For any integers $1 \le t \le l-1$ we have
\begin{eqnarray}
\lim_{{\scriptstyle L \rightarrow \infty} \atop
   {\scriptstyle L \equiv 0\, (n)}} q^{-\overline{E}(l\La_0, (t^L))}
F_{(t^L)}^{(l)}(q)
& = &
\sum_{\eta}
\frac{F^{(l-t)'}_{\eta}(q)\,
F^{(t)'}_{\eta}(q)}
{(q)_{\zeta_1}\cdots (q)_{\zeta_{n-1}} }, \label{eq:rklim1}\\
\eta &=& \left( (n-1)^{\zeta_{n-1}},\ldots ,(1)^{\zeta_{1}} \right),
\label{eq:rklim2}
\end{eqnarray}
where the ground state energy is
$\overline{E}(l\La_0, (t^L)) = \frac{tL(L-n)}{2n}$ and the sum runs over
the partitions $\eta$ satisfying $\eta_1 \le n-1$ as above
and $|\eta| \equiv 0$ mod $n$.
\end{proposition}
By Example \ref{ex:calB} (1) the same limit of $X^{(l)}_\mu(l\La_0)$ corresponds
the the choice ${\cal V} = V(t\La_0)$
in (\ref{vcal}).
Thus under the assumption (\ref{con:ffrk2}),
Proposition \ref{pro:rklim} and (\ref{xlimit}) implies
\begin{equation}
%\sum_i[V(t\La_0) \otimes V((l-t)\La_0): l\La_0 - i \delta]\, q^i
a^{V(t\La_0) \otimes V((l-t)\La_0)}_{l\La_0}(q)
= \sum_{\eta}
\frac{X^{(l-t)'}_{\eta}((l-t)\La_0)\,
X^{(t)'}_{\eta}(t\La_0)}
{(q)_{\zeta_1}\cdots (q)_{\zeta_{n-1}} }, \label{eq:rklim3}
\end{equation}
where the $\eta$-sum and the relation with $\zeta$ is
the same as in Proposition \ref{pro:rklim}.
This is an RSOS analogue of the spinon character formula
conjectured in \cite{NY2}.

\section{Discussion}\label{discussion}

Let us discuss further generalizations of
the results in Section \ref{sec:fflimit}.

\subsection{Fermionic string function for arbitrary $X^{(1)}_n$}
\label{discussionkuni}

The $q$-series
formulae in Proposition \ref{pro:kklim1},
Remark \ref{rem:fsg} and Proposition \ref{pro:kklim2}
are all originated in some limits of the unrestricted
1dsum $g_\mu(\la)$.
Let us discuss their possible extensions
to an arbitrary non-twisted affine Lie algebra $X^{(1)}_n$.
Our argument in this subsection will only concern
infinite series in $q$.
To seek their finite ($q$-polynomial) versions
as in Section \ref{ssec4.1}
is an important open problem.
So far only level 1 cases have been studied
for $B^{(1)}_n, C^{(1)}_n, D^{(1)}_n$ in \cite{DJKMO2} and
$A^{(2)}_{2n-1}, A^{(2)}_{2n}, D^{(2)}_{n+1}$ in \cite{KMOTU3}.
We hope to report higher level cases in near future.

Let $\alpha_1, \ldots, \alpha_n$ be the simple roots
of the classical subalgebra $X_n$.
For $1 \le a \le n$ set $t_a = 2/\vert \alpha_a \vert^2$
in the normalization $\vert \mbox{long root} \vert^2 = 2$.
Thus one has $\forall t_a = 1$ for simply laced
algebras and $t_a \in \{1,2, 3\}$ in general.
The root and coroot lattices are
denoted by $Q = \sum_{a=1}^n {\bf Z} \alpha_a$ and
$Q^{\vee} = \sum_{a=1}^n {\bf Z} t_a \alpha_a$, respectively.
Fix $s \in {\bf Z}_{\ge 1}$ arbitrarily.
For any positive integers $l_1, \cdots, l_s$ we set
$l = l_1 + \cdots + l_s$ and
$M^{(a)}_J = t_a(l_1 + \cdots + l_J)$ for $1 \le J \le s$ and
$1 \le a \le n$.
Consider the tensor product of the vacuum modules of $X^{(1)}_n$,
namely, ${\cal V} = \otimes_{J=1}^s V(l_J\La_0)$.
For its character we have
\begin{conjecture}\label{con:multikns}
Let $l, l_J, M^{(a)}_J\,  ( 1 \le J \le s)$ and ${\cal V}$ be as above.
For $\la \in Q$ we have
\begin{eqnarray}
%\sum_i (\mbox{dim }{\cal V}_{\la + l\La_0 - i\delta})\, q^i
q^{-\frac{(\la \vert \la)}{2l}}c^{\cal V}_\la(q)
& = &
\frac{1}{(q)_{\infty}^{n s}}
\sum_{\{ m\}^{\la}_{G_l} }
\frac{q^{C(\{ m \} )}}{\prod_{a=1}^{n}
\prod_{{\scriptstyle 1 \le i \le t_a l - 1} \atop
   {\scriptstyle i \notin \{M^{(a)}_1, \ldots, M^{(a)}_{s-1}\}}}
(q)_{m_{i}^{(a)}}}, \label{eq:multikns1}\\
C(\{ m \} ) & = &
\frac{1}{2} \sum_{G_l}
(\alpha_a \vert \alpha_b)\bigl(
\mbox{min}( t_b j, t_a k) - \frac{j k}{l} \bigr)
m^{(a)}_j m^{(b)}_k, \label{eq:multikns2}\\
G_l & = & \{ (a,j) \mid 1 \le a \le n, \, 1 \le j \le t_a l - 1\},
\label{eq:multikns3}
\end{eqnarray}
where the sum $\sum_{G_l}$ is taken over
$(a,j), (b,k) \in G_l$ and the sum $\sum_{\{ m\}^{\la}_{G_l} }$ is over
$\{ m^{(a)}_j \in {\bf Z}_{\ge 0} \mid (a,j) \in G_l \}$
satisfying the condition
$\sum_{(a,i) \in G_l} i\, m^{(a)}_i \alpha_a \equiv \la$ mod $l Q^{\vee}$.
\begin{eqnarray}
e^{-l\La_0} {\sl ch } {\cal V} & = &
\frac{1}{(q)_{\infty}^{n s}} \sum_{\{ m\}_{\overline{G}_l} }
\frac{q^{\overline{C}_l(\{ m \} )}}{\prod_{a=1}^{n}
\prod_{{\scriptstyle 1 \le i \le t_a l - 1} \atop
   {\scriptstyle i \notin \{M^{(a)}_1, \ldots, M^{(a)}_{s-1}\}}}
(q)_{m_{i}^{(a)}}}\ e^{\beta_l(\{m \})}, \label{eq:multikns4}\\
\overline{C}_l(\{ m \} ) & = &
\frac{1}{2} \sum_{\overline{G}_l}
(\alpha_a \vert \alpha_b) \mbox{min}( t_b j, t_a k)
m^{(a)}_j m^{(b)}_k, \label{eq:multikns5}\\
\beta_l(\{ m \}) & = &
\sum_{(a, i) \in \overline{G}_l} i \, m^{(a)}_i \alpha_a,
\label{eq:multikns6}\\
\overline{G}_l & = & G_l \sqcup_{a=1}^n \{ (a,t_a l) \},
\label{eq:multikns7}
\end{eqnarray}
where the sum $\sum_{\overline{G}_l}$ is taken over
$(a,j), (b,k) \in \overline{G}_l$ and
the sum $\sum_{\{ m\}_{\overline{G}_l} }$ is over
$\{ m^{(a)}_j \in {\bf Z}_{\ge 0} \mid (a,j) \in G_l \}$
and $\{ m^{(a)}_{t_a l} \in {\bf Z} \mid 1 \le a \le n \}$.
\end{conjecture}

For $X^{(1)}_n = A^{(1)}_{n-1}$,
(\ref{eq:multikns1})--(\ref{eq:multikns3}) reduce to
$\forall r_J = 0$ case of Proposition \ref{pro:kklim2},
and $s=1$ case of (\ref{eq:multikns4})--(\ref{eq:multikns7})
reduces to Remark \ref{rem:fsg}.
For $s = 1$ the conjecture
(\ref{eq:multikns1})--(\ref{eq:multikns3}) goes back to
\cite{KNS}.
In fact for any $X^{(1)}_n$ the two conjectures
(\ref{eq:multikns1})--(\ref{eq:multikns3}) and
(\ref{eq:multikns4})--(\ref{eq:multikns7}) are
equivalent.
Moreover they follow from the $s=1$ cases
by noting that
$e^{-l\La_0}{\sl ch } {\cal V} = \prod_{J=1}^s
e^{-l_J\La_0}{\sl ch } V(l_J\La_0)$.
Let us explain these facts more precisely.
To see the equivalence recall the decomposition of
characters in terms of theta functions (cf. eq.(12.7.12) in \cite{Ka}):
\begin{eqnarray}
e^{-l\La_0}{\sl ch } {\cal V} & = &
\sum_{\la \in Q/l Q^\vee} \theta_{\la, l}\ q^{-\frac{(\la \vert \la)}{2l}}
c^{\cal V}_\la(q),
%\sum_i (\mbox{dim }{\cal V}_{\la + l\La_0 - i\delta})\, q^i,
\label{chtheta}\\
\theta_{\la, l} & = & \sum_{\xi \equiv \la \bmod  l Q^\vee}
q^{\frac{(\xi \vert \xi)}{2l}} e^\xi.\label{thetafn}
\end{eqnarray}
Substitute (\ref{eq:multikns1}) and (\ref{thetafn}) into
(\ref{chtheta}).
The resulting double sum $\sum_{\la, \xi}$ is equivalent
to the single one $\sum_{\xi \in Q}$.
Moreover $\sum_{\xi \in Q} \sum_{\{m\}^\la_{G_l}}$
can further be replaced by $\sum_{\{m \}_{\overline{G}_l}}$
by identifying $\xi$ with $\beta_l(\{m\})$ in (\ref{eq:multikns6}).
On the other hand $\sum_{G_l}$ in
(\ref{eq:multikns2}) can be replaced with
$\sum_{\overline{G}_l}$ without causing any change.
Therefore
$C(\{m\}) = \overline{C}_l(\{m\}) -
\frac{(\beta_l(\{m\})\vert \beta_l(\{m\}))}{2l}$
holds between (\ref{eq:multikns2}) and (\ref{eq:multikns5}).
Combining these facts we obtain
(\ref{eq:multikns4})--(\ref{eq:multikns7}) from
(\ref{eq:multikns1})--(\ref{eq:multikns3}).
Obviously the converse of this argument is also valid.
To explain that Conjecture \ref{con:multikns} reduces to
$s = 1$, it suffices essentially to verify $s = 2$ case by induction.
Let us do so for (\ref{eq:multikns4})--(\ref{eq:multikns7}).
Taking the product of (\ref{eq:multikns4}) with
$(s,l) = (1,l_1)$ and $(1,l_2)$, we get the following
expression for
$e^{-(l_1+l_2)\La_0}{\sl ch } V(l_1\La_0) {\sl ch } V(l_2\La_0)$:

\begin{equation}
\frac{1}{(q)^{2n}_{\infty}}
\sum_{\{\tilde{m} \}_{\overline{G}_{l_1}}}
\sum_{\{\hat{m} \}_{\overline{G}_{l_2}}}
\frac{q^{\overline{C}_{l_1}(\{\tilde{m}\})+
\overline{C}_{l_2}(\{\hat{m}\})}}{
\prod_{(a,j) \in G_{l_1}}(q)_{\tilde{m}^{(a)}_j}
\prod_{(a,j) \in G_{l_2}}(q)_{\hat{m}^{(a)}_j}}
\, e^{\beta_{l_1}(\{\tilde{m}\})+\beta_{l_2}(\{\hat{m}\})}.
\label{eq:sistwo}
\end{equation}
For each color $1 \le a \le n$, replace here
as $\tilde{m}^{(a)}_{t_al_1} \rightarrow
\tilde{m}^{(a)}_{t_al_1}
+ \sum_{j=1}^{t_al_2} \hat{m}^{(a)}_j$ and introduce the
new variables
$\{ m^{(a)}_j \mid (a,j) \in \overline{G}_{l_1+l_2} \}$
by $m^{(a)}_j = \tilde{m}^{(a)}_j$ for $ 1 \le j \le t_al_1$,
$ = \hat{m}^{(a)}_{j-t_al_1}$ for $t_al_1 < j \le t_a(l_1+l_2)$.
Then it is straightforward to check
$\beta_{l_1}(\{\tilde{m}\})\vert_{\mbox{\scriptsize replacement}}
+\beta_{l_2}(\{\hat{m}\}) = \beta_{l_1 + l_2}(\{ m \})$ and
$\overline{C}_{l_1}(\{\tilde{m}\})\vert_{\mbox{\scriptsize replacement}}+
\overline{C}_{l_2}(\{\hat{m}\}) =
\overline{C}_{l_1+l_2}(\{m \})$.
Thus (\ref{eq:sistwo}) indeed yields
(\ref{eq:multikns4})--(\ref{eq:multikns7}) with $s = 2$
and $l = l_1 + l_2$.

\subsection{Fermionic form of $X^{(l)'}_\eta(\la)$
for $\la$ non vacuum type}\label{discussionhata}
In (\ref{con:ffrk2}), we conjectured a fermionic formula of
the restricted 1dsum $X^{(l)'}_\eta(\la)$
only for $\la = l\Lambda_0$.
Here we present a conjecture involving more general $\la$.
As the subsequent argument explains,
it is intimately related to the
spinon character formula (Conjecture \ref{spinonch1})
for general $\La$.

First we introduce a fermionic expression
$F^{(l,r)'}_{\eta,\mu}(q)$ ($0 \le r \le n-1$) for
partitions $\eta = ((n-1)^{\zeta_{n-1}},\ldots ,1^{\zeta_{1}})$
and $\mu$ satisfying $\mu_1\le l-1$ and
$|\eta| \equiv |\mu |+lr \bmod n$ by
\begin{eqnarray}
F^{(l,r)'}_{\eta,\mu}(q) & = &
\sum_{\{\tilde{m} \}} q^{\tilde{c}_l(\{ \tilde{m} \})}
\prod_{{\scriptstyle 1 \le a \le n-1} \atop
   {\scriptstyle 1 \le i \le l-1}}
\left[ \begin{array}{c} \tilde{p}^{(a)}_i +  \tilde{m}^{(a)}_i
 \\   \tilde{m}^{(a)}_i \end{array} \right],\label{eq:affrk7}\\
\tilde{c}_l(\{ \tilde{m}  \}) & = & \frac{1}{2}
\sum_{1 \le a, b \le n-1} C^{(n)}_{a b}
\sum_{1 \le j, k \le l} \mbox{min}(j, k)
\tilde{m}^{(a)}_j \tilde{m}^{(b)}_k
-\sum_{j=1}^{\mu'_1}\sum_{k=1}^l
\min (l-\mu_j,k)\tilde{m}^{(1)}_{k}  \nonumber\\
&& \quad
+ n(\mu) -\mu'_{1}(|\mu|+lr)
+ \frac{1}{n}(r+\mu'_1)(|\mu|+lr-|\eta|)-\frac{lr(r+1)}{2},\label{eq:affrk5}\\
\tilde{p}^{(a)}_i & = & C^{(l)^{-1}}_{1 i}\zeta_a
+ \delta_{a 1}\sum_{j=1}^{\mu'_1}C^{(l)^{-1}}_{l-\mu_j\, i}
- \sum _{b=1}^{n-1} C^{(n)}_{a b} \sum_{k=1}^{l-1}
C^{(l)^{-1}}_{i k} \tilde{m}^{(b)}_k,
\label{eq:affrk6}
\end{eqnarray}
Here the sum $\sum_{\{ \tilde{m} \}}$ runs over
$\{ \tilde{m}^{(a)}_i \in {\bf Z}_{\ge 0}
\mid 1 \le a \le n-1, \, 1 \le i \le l \}$
satisfying $\tilde{p}^{(a)}_i \ge 0$ for
$1 \le a \le n-1, \, 1 \le i \le l-1$, and
\begin{equation}\label{eq:affrk8}
\sum_{i=1}^l  i\,  \tilde{m}^{(a)}_i = \sum_{b=1}^{n-1}
C^{(n)^{-1}}_{a b} \zeta_b +l\mu'_1 - \frac{n-a}{n}|\mu|+\frac{alr}{n}
\quad \mbox{ for } \ 1 \le a \le n-1.
\end{equation}
Note that $F^{(l,r)'}_{\eta,\mu}(q)$ is a generalization of
$F^{(l)'}_{\eta}(q)$ in (\ref{eq:affrk1})--(\ref{eq:affrk4})
in that $F^{(l)'}_{\eta}(q) = F^{(l,0)'}_{\eta,\emptyset}(q)$.
Now our generalized conjecture reads
%%%%%%%%%%%%%%%%%%%%%%%%%%%%%%%%% Conj. of \sum q^E X %%%%%%%%%%%%
\begin{conjecture}\label{conj:r1dsumFF}
For any $0 \le r \le n-1$ and partitions $\eta, \mu$ satisfying
$\eta_1\le n-1,\;\mu_1\le l-1$ and
$|\eta| \equiv |\mu |+lr \mbox{ mod } n$, we have
\begin{equation}
\sum_{p \in {\cal H}(l\La_r,\mu)}q^{E(p)} X^{(l)'}_\eta(l\Lambda_r+af(\wt p))
= F^{(l,r)'}_{\eta,\mu}(q). \label{eq:hatahata}
\end{equation}
\end{conjecture}
\begin{remark}
When $\mu=(s)$ for some $0 \le s \le n-1$,
the set ${\cal H}(l\La_r,\mu)$ consists of
only one element $p$ as in Example \ref{ex:calB} (2), for which
$l\Lambda_r + af(\wt p) = (l-s)\Lambda_r+s\Lambda_{r+1}$ and $E(p)=0$
hold. Thus the above conjecture reduces to
\begin{equation}
X^{(l)'}_\eta((l-s)\Lambda_r+s\Lambda_{r+1}) =
F^{(l,r)'}_{\eta,(s)}(q)
\end{equation}
for any partition
$\eta$ satisfying $\eta_1\le n-1$ and $|\eta| \equiv lr+s \bmod n$.
\end{remark}

We explain briefly how the above conjecture emerged.
Let us consider the limit of the Kostka-Foulkes polynomial as
%%%%%%%%%%%%% Proposition : The limit of Kostka polynomial %%%%%%%%%%%%%%%%%%
\begin{proposition} \label{prop:klim2}
For finite partitions $\la$ and $\mu$
satisfying $l(\la)\le n-1,\;
\mu_1 \le l-1$ and $|\la|=|\mu|+lr$, define
$\tilde{\mu}=(l^L)\cup \mu, \;
\tilde{\la}=((\frac{|\tilde{\mu}|-|\la|}{n})^n) + \la$.
Then, we have
\begin{eqnarray}
\lim_{{\scriptstyle L \rightarrow \infty} \atop
   {\scriptstyle L \equiv r\, (n)}}
q^{-\overline{E}(l\Lambda_0,\tilde{\mu})}
K_{\tilde{\la} \, \tilde{\mu}}(q)
& = &
q^{-\overline{E}(l\Lambda_r,\mu)}
\sum_\eta
\frac{K_{\xi \eta}(q)\,
F^{(l,r)'}_{\eta,\mu}(q)}
{(q)_{\zeta_1}\cdots (q)_{\zeta_{n-1}} },\label{eq:klim5}\\
\xi = (n^{\frac{\vert \eta \vert - \vert \la \vert}{n}})
\cup \la',\quad
\eta &=& \left( (n-1)^{\zeta_{n-1}},\ldots ,1^{\zeta_{1}} \right),
\label{eq:klim7}\\
\nonumber
\end{eqnarray}
where the sum $\sum_\eta$ runs over the partitions
$\eta$ of the above form obeying
$|\eta| \equiv |\la|$ mod $n$.
\end{proposition}
This is obtained by taking the limit of (\ref{eq:ffk1}),
which is a very parallel calculation with Proposition \ref{pro:klim}.
Here we only mention that in its proof,
the ``minimum point'' $m^{(a)}_{i,0}$ should be replaced by
$m^{(a)}_{i, 0} = \frac{n-a}{n}(L-r)\delta_{i l}$
and the relation
$\overline{E}(l\Lambda_0,\tilde{\mu})=\overline{E}(l\Lambda_r,\mu)
+\frac{l}{2n}(L^2-r^2) + (\frac{|\mu|}{n}-\frac{l}{2})(L-r)$
is used.  See (\ref{gsenergy}).

{}From (\ref{xbarbyk}), (\ref{xbarlimit}) and (\ref{xbarlimit2}),
we know that
\begin{equation}
\mbox{LHS of (\ref{eq:klim5})} =
\sum_{p \in {\cal H}(l\La_r, \mu)}
q^{E(p) - \overline{E}(l\Lambda_r,\mu)}
b^{V(l\La_r + af(\wts p))}_\la(q). \nonumber
\end{equation}
Admitting Conjecture \ref{spinonch1}, we rewrite this as
\begin{equation}
q^{- \overline{E}(l\Lambda_r,\mu)}
\sum_{p \in {\cal H}(l\La_r, \mu)}
q^{E(p)} \sum_\eta
\frac{X'_{\eta}(\la)\,
X^{(l)'}_{\eta}(l\Lambda_r+af(\wt p))}
{(q)_{\zeta_1}\cdots (q)_{\zeta_{n-1}} },\label{eq:klim4} \\
\end{equation}
where the sum $\sum_\eta$ and $\zeta$ are specified
as in Proposition \ref{prop:klim2}.
{}From (\ref{xbarbyk}) we have  $X'_\eta(\la) = K_{\xi \eta}(q)$
for $\xi$ given in (\ref{eq:klim7}).
Thus by comparing (\ref{eq:klim4}) with (\ref{eq:klim5}),
we arrive at Conjecture \ref{conj:r1dsumFF} naturally.

\begin{example}
Let $n=3$, $r=1$, $l=3$, $\eta =(2,2,1,1)$, $\mu =(2,1)$. Then,
${\cal H}(l\Lambda_r,\mu)=\{ 22 \otimes 2,22\otimes 3\}$, where
the element $(x_1,x_2,x_3)$ of each crystal $B_{(k)}$ is denoted by
$1^{x_1}2^{x_2}3^{x_3}$.
The LHS of (\ref{eq:hatahata}) is calculated from the data
\begin{center}
\begin{tabular}{|c|c|c|c|c|}\hline
$p \in {\cal H}(l\Lambda_1,\mu)$ & $E(p)$ & $3 \Lambda_1 + \wt p$ &
$p_0 \in {\cal P}^{(3)'}_\eta(3\Lambda_1 + \wt p)$ & $E(p_0)$
\\ \hline \hline
$22 \otimes 2$ & $1$ & $3 \Lambda_2$ & $12 \otimes 12 \otimes \otimes 1 \otimes 2$
& $-1$ \\ \hline
$22 \otimes 3$ & $0$ & $3 \Lambda_0+\Lambda_1+\Lambda_2$ &
$12 \otimes 12 \otimes \otimes 1 \otimes 3$ & $-1$ \\
&&& $12 \otimes 12 \otimes \otimes 3 \otimes 1$ & $-2$ \\
&&& $12 \otimes 13 \otimes \otimes 1 \otimes 2$ & $-2$ \\
&&& $12 \otimes 13 \otimes \otimes 2 \otimes 1$ & $-3$ \\ \hline
\end{tabular}
\end{center}
as
\begin{equation}
\makebox{LHS} = q^1q+q^0(q+2q^2+q^3) = q+3q^2+q^3.
\end{equation}
On the other hand, $F^{(3,1)'}_{\eta ,\mu}(q)$ defined
by (\ref{eq:affrk7})-(\ref{eq:affrk8}) is calculated from the data below,
which illustrates
all the admissible values of $\{ \tilde{m}_{j}^{(a)}\}$ for $\sum_{\{\tilde{m}\} }$
in (\ref{eq:affrk7}) and the corresponding $\{ \tilde{p}_{j}^{(a)}\}$, etc,
\begin{center}
{\small
\begin{tabular}{|cc|cc|c|cc|cc|c||c|} \hline
$\tilde{m}_{1}^{(1)}$ & $\tilde{p}_{1}^{(1)}$ & $\tilde{m}_{2}^{(1)}$ & $\tilde{p}_{2}^{(1)}$ & $\tilde{m}_{3}^{(1)}$ &
$\tilde{m}_{1}^{(2)}$ & $\tilde{p}_{1}^{(2)}$ & $\tilde{m}_{2}^{(2)}$ & $\tilde{p}_{2}^{(2)}$ & $\tilde{m}_{3}^{(2)}$ &
$\tilde{c}_3(\{ \tilde{m} \})$ \\ \hline \hline
$0$ & $2$ & $2$ & $0$ & $1$ & $1$ & $0$ & $1$ & $0$ & $2$ & $2$ \\ \hline
$1$ & $1$ & $0$ & $1$ & $2$ & $0$ & $2$ & $0$ & $1$ & $3$ & $1$ \\ \hline
$2$ & $0$ & $1$ & $0$ & $1$ & $1$ & $1$ & $1$ & $0$ & $2$ & $2$ \\ \hline
\end{tabular}
}
\end{center}
as
\begin{equation}
F^{(3,1)'}_{\eta ,\mu}(q)=q^2+q(1+q)+q^2(1+q)=q+3q^2+q^3.
\end{equation}
\end{example}

\vskip0.5cm\noindent
{\bf Acknowledgements}.
A.~K., A.~N.~K. and M.~O. thank the organizers of
{\em Workshop on Algebraic Combinatorics},  June 9--20, 1997
held at CRM, University of Montreal
for hospitality, and A.~Lascoux, B.~Leclerc and J.-Y.~Thibon for
discussions and correspondence on the program ``kostka40.c''.
A.~N.~K. wishes to thank L.~Vinet for hospitality at the CRM,
University of Montreal, where part of this work was done.
A.~K. thanks the organizers of {\em International workshop on
statistical mechanics and integrable systems},
July 20 -- August 8, 1997 held in Coolangatta and Canberra,
where part of this work was presented.
He thanks P.~Bouwknegt, R.~Kedem, B.~M.~McCoy and A.~Schilling
for discussions and comments.

\end{document}